\documentclass[12pt,a4paper]{article}
\usepackage{latexsym}
\usepackage{amsmath,amssymb}
\usepackage{color}
\usepackage{epsfig,graphicx,float}
\usepackage{authblk}

\newcommand{\be}{\begin{equation}}
\newcommand{\ee}{\end{equation}}

\newtheorem{proposition}{Proposition}[section]
\newtheorem{theorem}[proposition]{Theorem}

\newtheorem{lemma}[proposition]{Lemma}
\newtheorem{definition}[proposition]{Definition}
\newtheorem{remark}[proposition]{Remark}

\newcommand{\R}{I\!\!R}
\newcommand{\N}{I\!\!N}

\newcounter{secnum}

\title{Inverse source problem in a forced network}

\author[1,2]{J. G. Caputo \thanks{caputo@insa-rouen.fr}}
\author[1]{A. Hamdi \thanks{adel.hamdi@insa-rouen.fr}}
\author[1]{A. Knippel \thanks{arnaud.knippel@insa-rouen.fr}}

\affil[1]{Laboratoire de Math\'ematiques, INSA de Rouen Normandie\\ 
76801 Saint-Etienne du Rouvray, France.}

\affil[2]{Laboratoire de Math{\'e}matiques Rapha{\"e}l Salem, Universit{\'e}
de Rouen Normandie, 76801 Saint-Etienne-du-Rouvray,  France.}

\date{\ }

\begin{document}
\maketitle

\begin{abstract} 
We address the nonlinear inverse source problem of identifying a 
time-dependent source occurring in one
node of a network governed by a wave equation. We prove that 
time records of the associated state taken at a strategic set 
of two nodes yield
uniqueness of the two unknown elements: the source position and the 
emitted signal. We develop a non-iterative identification method that 
localizes the source node by solving
a set of well posed linear systems. Once the source node is localized, 
we identify the emitted signal using a deconvolution problem 
or a Fourier expansion.
Numerical experiments on a $5$ node graph confirm the effectiveness of the approach.
\end{abstract}

\section{Introduction}

Networks play a crucial role in the transmission of numerous 
quantities. These can be miscible like fluids \cite{Maas} 
and electric power \cite{Chow}, or non miscible like 
information packets \cite{Tsitsiklis} in telecommunications. The 
mathematical structure describing a network is a graph composed 
of nodes connected by edges. The graph can be directed as for
traffic flow or nerve propagation \cite{scott} or undirected 
as for miscible flows; in these two situations the equations
differ.

Miscible flows satisfy discrete conservation laws; for electrical 
circuits, these are the well-known Kirchoff and Ohm's laws. When 
dissipation is absent from the edges, miscible flows obey 
a "graph wave equation" where the graph Laplacian matrix replaces the 
ordinary Laplacian \cite{Maas,mohar91}. This is a very general
model, essentially Newton's equation, linking acceleration to
the forces applied on the system. It is valid for acoustic waves 
in a fluid or normal modes in a molecule \cite{landau}. It is a 
linear limit of the so-called Kuramoto 
model for an electrical network \cite{timme}. In that context, the
variable is the phase of a generator or load, located at each node of
the network. Note that one could introduce damping at each node
of the system
leading to a so-called reaction-diffusion type dynamics \cite{scott}.

For such a miscible flow network, the forward problem consists in finding 
the state of the transmitted entity in each node of the network, given 
the physical parameters and the initial state. For many applications, the 
inverse problem is important. There, we aim to identify some unknown 
parameters affecting the state of the system using measurements 
at some observation nodes of the network. For example, we would like to detect 
the presence of a failing node on the network -typically an unsynchronized
generator in the electrical grid context- and possibly 
localize its position to take appropriate action. 
Many authors analyzed source problems in networks; 
two examples are the articles by Rodriguez and Sahah 
\cite{Rodriguez,Sahah} who 
-in an effort to contain an epidemic or contamination-
analyzed the influence of the propagation network and of
the unknown source, assuming known the state in all the nodes.

In practice, it is usually infeasible to record the state of the system at all the nodes of 
the network and recent studies address this constraint. That on a 
used network only a few number of its nodes can be observed. For example, Pinto et al \cite{Pedro} propose a probabilistic approach to estimate the location of a node source using measurements collected at sparsely-placed observation nodes in the network. The authors conclude that such a sparse deployment of sensors provides an effective alternative to recording the state at all nodes of the network. However, the position of the observation nodes strongly affects the performance of the approach and thus optimal strategies for placing the sensors need to be found. In \cite{Kai}, the authors developed the so-called Short-Fat Tree algorithm to identify a node diffusion source in networks; they proved limit results on the asymptotic probability of their approach in identifying the source node with respect to the infection duration. In the power systems context, Nudell and Chakrabortty \cite{Nudell} propose a graph theoretic algorithm to find forced oscillation  inputs. Minimal information is assumed to be known in order to estimate the transfer function.

In the present study, we consider a dissipation less graph wave
equation modeling a general miscible flow, forced by an unknown
signal generated at an unknown node. We develop a constructive deterministic approach 
based on time records of the transmitted entity in a strategic set 
of two nodes of the network. We show that, with this method,  we localize 
the node source and identify its unknown emitted signal in a 
unique manner. We introduce an adjoint problem formed by convoluting the 
equation with an appropriate boundary value problem, yielding an 
over-determined system of $N$ equations and $N-2$ unknowns for 
a $N$ node network. The right hand side of this system contains 
the unknown final state $X(T)$ of the network. A first step is 
to estimate $X(T)$; for this it is crucial that the observation nodes form 
a strategic set. Once $X(T)$ is found, we show how to solve the 
over-determined linear system and find the node source. Graph theory
arguments provide guidelines on how to choose the observation nodes
for a given network.
Then we give two methods to estimate the unknown signal. \\
The paper is organized as follows. Section two introduces the model and the mathematical framework. Section 3 presents the adjoint problem and the over-determined linear system. We establish the identifiability of the source location and signal in section 4 and develop an identification method in section 5. 
Numerical calculations on a 5 node graph are shown in section 6 
and we conclude in section 7.

\section{Mathematical Modelling and technical results}

Consider a network defined over N nodes labelled $k=1,\dots,N$ 
on which propagates a miscible entity. The state of this entity 
at time $t \in (0,T)$ in the different nodes of the network 
is given by the following vector
\begin{eqnarray}
X(t)=\big(x_1(t),\dots,x_N(t)\big)^{\top} ,
\label{state}
\end{eqnarray}      
where $T>0$ is a final monitoring time. Neglecting friction in the nodes and 
the edges, the evolution of the state $X$ is governed by 
the so-called graph wave equation, see \cite{Maas}\cite{cks13}: 
\begin{eqnarray}
\left\{
\begin{array}{lll}
 {\ddot X}(t) - \Delta X(t) = \lambda(t) S \quad \mbox{in} \; (0,T)\\
X(0)=a\in  \R^N \quad \mbox{and} \quad \dot X(0)=b\in \R^N,
\end{array}
\right.
\label{weq}
\end{eqnarray}
where $a$ and $b$ represent, respectively, the initial state and 
velocity in all nodes of the network, $\Delta$ is the 
$N\times N$ graph Laplacian matrix \cite{bls07} and $\lambda S$ 
designates a node source of the network located at $s$ and forcing 
the state of the transmitted entity by emitting a time-dependent 
signal $\lambda$. 
We assume $\lambda$ and $S$ to be elements of the admissible set:
\begin{eqnarray*}
\begin{array}{lll}
{\cal{A}}:=\big\{\lambda\in L^2(0,T) \;\; \mbox{that satisfies:} \;\; \exists T^0\in(0,T)/ \; \lambda(t)=0, \; \forall \; t\in(T^0,T) \quad \mbox{and}\\
\qquad S=\big(s_1,\dots,s_N\big)^{\top} \;\; \mbox{where for} \;k=1,\dots,N,\;\; s_k\in\{0,1\} \quad \mbox{and} \quad \displaystyle \sum_{k=1}^N s_k=1\big\}
\end{array}
\label{Ad}
\end{eqnarray*}

Since the graph wave equation (\ref{weq}) is a system of linear second 
order differential equations it follows that for all given two admissible 
elements $\lambda$ and $S$ defining a node source $\lambda S$, the system 
(\ref{weq}) admits a unique solution $X=(x_1,\dots,x_N)^{\top}$ such 
that $x_k\in H^2(0,T)$ for $k=1,\dots,N$. That allows to define the 
following observation operator:
\begin{eqnarray}
M\big[\lambda,S\big]:=\big\{x_i(t),x_{j}(t); \;\mbox{for} \; 0<t<T\big\} ,
\label{observations}
\end{eqnarray}
where $i\ne j \in\{1,\dots,N\}$ are two distinct nodes of the network. This 
is the so-called {\it{forward problem}}.

The {\it{inverse problem}} with which we are concerned here is: given time 
records $d_i$ and $d_j$ in $(0,T)$ of the two local state 
components $x_i$ and $x_j$ taken at the two distinct 
nodes $i$ and $j$ of the network, determine the two unknown 
elements $\lambda$ and $S$ of ${\cal{A}}$ that yield
\begin{eqnarray}
M\big[\lambda,S\big]=\big\{d_i(t),d_j(t); \;\mbox{for} \; 0<t<T\big\}
\label{records}
\end{eqnarray}

Besides, the graph Laplacian matrix $\Delta$ is a $N\times N$ 
real symmetric and negative semi-definite matrix, for more details 
see \cite{bls07}. In the remainder, we assume the eigenvalues 
of $\Delta$ to be $N$ distinct nonpositive real numbers ordered and 
denoted as follows:
\begin{eqnarray}
0=- \omega_1^2 > - \omega_2^2 > \dots > -\omega_N^2 .
\label{eigenvalues_Delta}
\end{eqnarray}
The set of eigenvectors $\big\{v^1,\dots,v^N\big\}$ associated to these
eigenvalues is an orthogonal base of $\R^N$ and can be chosen
orthonormal. Thus, the state $X$ solution of the system (\ref{weq}) 
can be written in $(0,T)$ as follows:
\begin{eqnarray}
X(t)=\displaystyle\sum_{n=1}^{N}y_n(t)v^n \qquad \qquad \mbox{where} \quad y_n(t)=\langle X(t),v^n\rangle \qquad \mbox{for} \; n=1,\dots,N ,
\label{X_in_vn}
\end{eqnarray}
where $\langle,\rangle$ designates the inner product in $\R^N$.

We now introduce the definition of a {\it{strategic}} set of nodes.
\begin{definition} 
A set of $\{k_1,\dots,k_{\ell}\}$ nodes is called {\it{strategic}} if for 
all vector $v^n$, where $n=1,\dots,N$ there exists at least one element 
$k\in \{k_1,\dots,k_{\ell}\}$ such that $v^n_k \ne 0$.
\label{strategic} 
\end{definition}


Therefore, according to Definition \ref{strategic}, a single node 
$i\in\{1,\dots,N\}$ is called {\it{strategic}} if $v^n_i\ne 0$ for 
all $n=1,\dots,N$. This is the opposite of a "soft node" that was introduced
in our previous work \cite{cks13}. Also, for spatially continuous systems
Jai and Pritchard stressed the importance of these strategic sets \cite{Jai}.

A set of two distinct nodes $\{i,j\}$ is {\it{strategic}} if 
for all $n=1,\dots,N$, it holds $|v^n_i| + |v^n_j|\ne 0$. Hence, 
if $i$ is a {\it{strategic}} node then, all set of nodes containing 
$i$ is also {\it{strategic}}. This leads to the following technical result:


\begin{lemma}
Let $T^*\in (T^0,T)$, $X$ be fulfilling the following system:
\begin{eqnarray}
\left\{
\begin{array}{lll}
{\ddot X}(t) - \Delta X(t)=0 \quad \mbox{in} \; (T^*,T)\\
X(T) \in \R^N \quad \mbox{and} \quad \dot X(T)\in \R^N
\end{array}
\right. ,
\label{Lemma_1}
\end{eqnarray}
and $\{k_1,\dots,k_{\ell}\}$ be a {\it{strategic}} set of nodes. We have 
\begin{eqnarray}
x_{k_1}(t)=\dots=x_{k_{\ell}}(t)=0, \;\; \forall \;t\in (T^*,T) \implies  X(T)=\dot X(T)=\vec 0
\label{Result_lemma}
\end{eqnarray}
\label{Null_State}
\end{lemma}

{\bf{Proof.}} See the appendix.

\section{Adjoint problem}

To solve the inverse problem for (\ref{weq}) , we
introduce an adjoint problem.
For that, we consider the Sturm-Liouville eigenvalue problem:
\begin{eqnarray}
\left\{
\begin{array}{lll}
-\ddot\varphi_m(t)=\mu_m \varphi_m(t) \quad \mbox{in} \; (0,T)\\
\varphi_m(0)=\varphi_m(T)=0 .
\end{array}
\right.
\label{Sturm}
\end{eqnarray}
Then, the normalized eigenfunctions $\varphi_m$ and their associated eigenvalues $\mu_m$ are defined for all $m\in \N^*$ by
\begin{eqnarray}
\varphi_m(t)=\displaystyle\sqrt{\frac{2}{T}}\displaystyle\sin\big(\displaystyle\frac{m\pi}{T}t\big) \qquad \mbox{and} \qquad \mu_m=\Big(\displaystyle\frac{m\pi}{T}\Big)^2>0
\label{eigen}
\end{eqnarray}
The set $\{\varphi_m\}$ forms a complete orthonormal 
family of $L^2(0,T)$ for the standard inner product in $L^2(0,T)$ i.e., 
$\langle f,g\rangle_{L^2(0,T)}=\int_0^T f(t)g(t)dt$.

To obtain the adjoint problem, we multiply each equation of the 
system (\ref{weq}) by $\varphi_m$ and integrate by parts 
over $(0,T)$, to get the following linear system:
\begin{eqnarray}
-\big(\Delta + \mu_m {\bf{I}}\big) \bar X_m=\lambda_m S + P_m ,
\label{bar_Xm}
\end{eqnarray}
where {\bf{I}} is the $N\times N$ identity matrix and 
\begin{eqnarray}
\begin{array}{lll}
\bar X_m=\big(\bar x_1 ,\dots,\bar x_N\big)^{\top} \quad \mbox{with} \quad \bar x_k=\langle x_k,\varphi_m\rangle_{L^2(0,T)}\\
P_m=\dot\varphi_m(T) X(T) - \dot\varphi_m(0) X(0) \qquad \mbox{and} \qquad \lambda_m=\langle \lambda ,\varphi_m\rangle_{L^2(0,T)}
\end{array}
\label{Notations}
\end{eqnarray}

Notice that since the final state $X(T)$ of the system (\ref{weq}) is 
unknown, the first term defining the vector $P_m$ introduced 
in (\ref{Notations}) is also unknown. Estimating $X(T)$ will 
then be the first step in the identification process.

Furthermore, the choice of 
the boundary conditions on $\varphi_m$ solution of the auxiliary 
problem (\ref{Sturm}) allowed to eliminate from the linear 
system (\ref{bar_Xm}) the unknown term $\dot X(T)$. If we had chosen
other boundary conditions for $\varphi_m$ we would always have 
one unknown term among $X(T)$ and $\dot X(T)$.

To simplify our presentation we introduce for all $m\in \N^*$ the notation
\begin{eqnarray}
A_m:=\Delta + \mu_m {\bf{I}}
\label{mat_A}
\end{eqnarray}
Furthermore, according to the observation operator $M[\lambda,S]$ in (\ref{observations}) it follows that the two components $\bar x_i=\langle x_i,\varphi_m\rangle_{L^2(0,T)}$ and $\bar x_j=\langle x_j,\varphi_m\rangle_{L^2(0,T)}$ of the vector $\bar X_m$ are known. Thus, (\ref{bar_Xm}) is reduced to the 
following $N\times (N-2)$ linear system:
\begin{eqnarray}
-A_m^{i,j} \bar X^{i,j}_m=P^{i,j}_m +\lambda_m S ,
\label{adj1}
\end{eqnarray}
where $A_m^{i,j}$ is the $N\times(N-2)$ matrix obtained by removing the two columns $i$ and $j$ from the $N\times N$ matrix $A_m$ introduced in (\ref{mat_A}), $\bar X^{i,j}_m \in \R^{N-2}$ is the unknown vector defined by removing the two known components $\bar x_i$ and $\bar x_j$ from the vector $\bar X_m$ in (\ref{Notations}) and
\begin{eqnarray}
P^{i,j}_m=\dot\varphi_m(T) X(T) - \dot\varphi_m(0) X(0) + A_m(:,i) \bar x_i + A_m(:,j) \bar x_j ,
\label{Pm_ij}
\end{eqnarray}
where $A_m(:,i)$ is the $i^{th}$ column vector of the 
matrix $A_m$ in (\ref{mat_A}). {At this stage, $\lambda_m, ~S$ and
$P^{i,j}_m$ are unknown, they need to be determined to solve
the linear system \eqref{adj1}. }

\section{Identifiability}

In this section, we prove that -under some reasonable assumptions-
the observation operator $M[\lambda,S]$ introduced in (\ref{observations}) 
determines uniquely the two unknown elements $\lambda$ and $S$ of 
the admissible set ${\cal{A}}$ defining the source $\lambda S$ 
occurring in the system (\ref{weq}). 

\subsection{Estimation of $X(T)$ }

As indicated above, the first step is to determine the unknown 
final state $X(T)$ defining the vector $P_m$ in 
(\ref{bar_Xm})-(\ref{Notations}) {from the measurements
$d_k(t)$ at the nodes $k$ of a strategic set.} We prove the following result:
\begin{theorem} Provided $\lambda$ and $S$ are two elements of the admissible set ${\cal{A}}$, time records in $(0,T)$ of $X$ the solution of (\ref{weq}) taken in a strategic set of nodes identify uniquely $X(T)$.   
\label{LS}
\end{theorem}

{\bf{Proof.}} For $n=1,2$, let $\lambda^{(n)}$ and $S^{(n)}$ be two elements of the admissible set ${\cal{A}}$ and $X^{(n)}$ be the solution of the system (\ref{weq}) with the node source $\lambda^{(n)}S^{(n)}$ instead of $\lambda S$. Then, the variable $X=X^{(2)}-X^{(1)}$ solves the system:

\begin{eqnarray}
\left\{
\begin{array}{lll}
{\ddot X}(t) - \Delta X(t)=\lambda^{(2)}(t)S^{(2)}-\lambda^{(1)}(t)S^{(1)} \quad \mbox{in} \; (0,T)\\
X(0)=\dot X(T)=0\in \R^N
\end{array}
\right.
\label{Th_1}
\end{eqnarray}
Since $\lambda^{(2)}$ and $\lambda^{(1)}$ are two admissible time-dependent source intensities, it follows that $\lambda^{(2)}=\lambda^{(1)}=0$ in $(T^0,T)$. Thus, the variable $X$ satisfies the system (\ref{Lemma_1}) for all $T^*\in(T^0,T)$. Moreover, if the local time records of $X^{(2)}$ and of $X^{(1)}$ 
taken in a strategic set of nodes coincide in $(0,T)$ then the 
values of the state $X$ in all nodes of a strategic set is null in $(0,T)$. Therefore, from applying Lemma \ref{Null_State} we conclude that $X(T)=\vec 0$ which means $X^{(2)}(T)=X^{(1)}(T)$. \hspace{6cm} $\blacksquare$

Then, using the solution of (\ref{Lemma_1}) given in (\ref{Express_X}) it follows that $X(t)=\big(x_1(t),\dots,x_N(t)\big)^{\top}$ is redefined for all $t\in (T^*,T)$ by
\begin{eqnarray}
\begin{array}{llll}
X(t)=\Big(y_1(T) + (t-T)\dot y_1(T)\Big)v^1  + \displaystyle \sum_{n=2}^{N} \Big(y_n(T)\cos\big(\omega_n(t-T)\big) + \displaystyle \frac{\dot y_n(T)}{\omega_n}\sin\big(\omega_n(t-T)\big)\Big)v^n
\end{array}
\label{Express_X_ident}
\end{eqnarray}

According to Theorem \ref{LS}, we determine the coefficients $y_n(T)$ and $\dot y_n(T)$ for $n=1,\dots,N$ defining in (\ref{Express_X_ident}) the state $X$ in $(T^*,T)$ by solving the minimization problem:
\begin{eqnarray}
\displaystyle \min_{Y\in \R^{2N}} \frac{1}{2} \displaystyle\sum_{k\in {\cal{E}}}\big\|x_k-d_k\big\|_{L^2(T^*,T)}^2 
\label{Least_squares}
\end{eqnarray} 
{where $x_k(t)$ is the $k$th component of $X(t)$,}
$Y=\big(y_1(T),\dots,y_N(T),\dot y_1(T),\dots,\dot y_N(T)\big)^{\top}$
and ${\cal{E}}$ is a strategic set of nodes.

\begin{remark} \label{remark_ls}
From the proof of Theorem \ref{LS}, {it follows that a single
strategic node is enough to identify uniquely $X(T)$ via
the minimization problem (\ref{Least_squares}).}
\end{remark}

\subsection{Identifiability of the source}

This result is given by the following identifiability theorem: 
\begin{theorem} 
Let $m\in \N^*$ and $A_m$ be the $N\times N$ matrix introduced in 
(\ref{mat_A}). Let $i,j$ be two nodes. If
\begin{enumerate}
\item the set $\{i,j\}$ is strategic,
\item $\langle \lambda,\varphi_m\rangle_{L^2(0,T)}\ne 0$,
\item all $(N-2)\times(N-2)$ matrices defined by removing from $A_m$ its $i^{th}, j^{th}$ columns and any $2$ rows are invertible,
\end{enumerate}
then, the observation operator $M[\lambda,S]$ in (\ref{observations}) 
identifies uniquely the two unknown elements $S$ and $\lambda$ of the 
admissible set ${\cal{A}}$ defining the source $\lambda S$ occurring in the system (\ref{weq}).  
\label{identifiability} 
\end{theorem}

{\bf{Proof.}} For $n=1,2$, let $\lambda^{(n)}$ and $S^{(n)}$ be 
two elements of ${\cal{A}}$ and $X^{(n)}$ be the solution of (\ref{weq}) with the source $\lambda^{(n)}S^{(n)}$ instead of $\lambda S$. We denote $X=(x_1,\dots,x_N)^{\top}$ the variable defined by $X=X^{(2)}-X^{(1)}$. Then, $X$ solves the system (\ref{Th_1}) and in view of (\ref{observations}), we have
\begin{eqnarray}
M\big[\lambda^{(2)},S^{(2)}\big]=M\big[\lambda^{(1)},S^{(1)}\big] \quad \implies \quad x_i(t)=x_j(t)=0 \;\forall t\in (0,T)
\label{records_ident}
\end{eqnarray}

The right hand side of the implication in (\ref{records_ident}) leads to 
the two following results: \\
{\bf{1.}} Applying Lemma \ref{Null_State}, we get $X(T)=\dot X(T)=\vec 0$. \\
{\bf{2.}} From (\ref{Notations}), we find $\bar x_i=\bar x_j=0$. \\
Therefore, the adjoint system associated to $X$ is
\begin{eqnarray}
-A_m^{i,j} \bar X^{i,j}_m=\lambda_m^{(2)}S^{(2)} - \lambda_m^{(1)}S^{(1)} ,
\label{bar_Xm_ident}
\end{eqnarray}
where $\lambda^{(n)}_m=\langle \lambda^{(n)} ,
\varphi_m\rangle_{L^2(0,T)}\ne 0$ for $n=1,2$. 
The linear system (\ref{bar_Xm_ident}) contains $N$ equations with 
$(N-2)$ unknowns. Furthermore, since $S^{(n=1,2)}$ belong to ${\cal{A}}$, 
the right hand side in (\ref{bar_Xm_ident}) contains at most two 
non-null components. Therefore, by taking out two equations containing 
the two eventual non-null right hand side terms, we end up with an 
homogeneous $(N-2)\times (N-2)$ linear system. Then, the
$(N-2) \times (N-2)$ matrix is invertible from the third
hypothesis of Theorem \ref{identifiability}, so that
$\bar X^{i,j}_m=\vec 0$ and thus,
\begin{eqnarray}
\lambda_m^{(2)}S^{(2)}=\lambda_m^{(1)}S^{(1)} \qquad \implies \quad S^{(2)}=S^{(1)} \quad\mbox{and} \quad \lambda_m^{(2)}=\lambda_m^{(1)}
\label{identifiab_S}
\end{eqnarray}

We set $S^{(2)}=S^{(1)}=S$ in the system (\ref{Th_1}). To prove that we have also $\lambda^{(2)}=\lambda^{(1)}$ in $(0,T)$, we expand the solution 
$X$ of (\ref{Th_1}) in the orthonormal family $\{v^1,\dots,v^N\}$ as 
in (\ref{X_in_vn}) i.e., 
$$X(t)=\sum_{n=1}^N y_n(t) v^n . $$ 
Then, $y_{n=1,\dots,N}$ solve 
\begin{eqnarray}
\left\{
\begin{array}{lll}
\ddot y_n(t) + \omega^2_n y_n(t)=\big(\lambda^{(2)}(t) - \lambda^{(1)}(t)\big) \langle S,v^{n}\rangle \quad \mbox{in} \; (0,T)\\
y_n(0)=\dot y_n(0)=0
\end{array}
\right.
\label{Ident_Lam}
\end{eqnarray}

Besides, we determine the two following fundamental solutions:
\begin{eqnarray}
\begin{array}{llll}
y_1^0(t)={\cal{H}}(t) t \qquad\qquad\qquad\quad \mbox{solves} \quad \ddot y_1^0(t)=\delta(t)\\
y_n^0(t)=\displaystyle\frac{1}{\omega_n}{\cal{H}}(t)\sin(\omega_n t) \qquad\; \mbox{solves} \quad \ddot y^0_n(t) + \omega_n^2 y_n^0(t)=\delta(t), \quad \forall n=2,\dots,N
\end{array}
\label{Fundamental}
\end{eqnarray}
where ${\cal{H}}$ is the Heaviside function and $\delta$ is the Dirac mass. Thus, (\ref{Fundamental}) gives the solution of (\ref{Ident_Lam}) for all $n=1,\dots,N$ and then, $X$ the solution of (\ref{Th_1}) for $S^{(2)}=S^{(1)}=S$ as follows:

\begin{eqnarray}
\begin{array}{lll}
X(t)=\langle S,v^{1}\rangle \displaystyle \int_0^t \big(\lambda^{(2)}(s) - \lambda^{(1)}(s)\big)(t-s)ds\; v^1\\
\qquad\qquad + \displaystyle \sum_{n=2}^{N} \Big(\frac{\langle S,v^{n}\rangle}{\omega_n}\displaystyle \int_0^t \big(\lambda^{(2)}(s) - \lambda^{(1)}(s)\big)\sin(\omega_n(t-s))ds\Big) v^n  \qquad \mbox{in} \; (0,T)
\end{array}
\label{X_part}
\end{eqnarray}
 
As from (\ref{records_ident}) we have $x_{k=i,j}(t)=0$ for all $t\in (0,T)$, it follows in view of (\ref{X_part}) that

\begin{eqnarray}
\begin{array}{lll}
\displaystyle\int_0^t \big(\lambda^{(2)}(s) - \lambda^{(1)}(s)\big) \Phi_{k=i,j}(t-s)ds=0, \qquad \forall \; t\in(0,T)\\
\mbox{where:} \qquad \Phi_k(t)=\langle S,v^{1}\rangle t  v^1_k + \displaystyle\sum_{n=2}^{N} \frac{\langle S,v^{n}\rangle}{\omega_n} \sin(\omega_nt) v_k^n
\end{array}
\label{convolution_identifiab}
\end{eqnarray}

Suppose that $\Phi_{k=i,j}= 0$ in $(t_k,\bar t_k)$ where $0<t_k<\bar t_k<T$. Then, using similar techniques as employed in (\ref{Express_X})-(\ref{dot_X_initial}) it comes from (\ref{convolution_identifiab}) that $\langle S,v^{n}\rangle v^n_{k=i,j}=0$ for all $n=1,\dots,N$. That's absurd since the set $\{i,j\}$ is strategic. Hence, it follows that for $k=i$ and/or $k=j$ we have $\Phi_k\ne 0$ a.e. in (0,T). Therefore, from (\ref{convolution_identifiab}) and according to Titchmarsh Theorem on convolution of $L^1$ functions it comes that $\lambda^{(2)}=\lambda^{(1)}$ in $(0,T)$. \hspace{1cm}$\blacksquare$


\subsection{Appropriate observation nodes}

The third condition of Theorem \ref{identifiability} i.e.
that all $(N-2)\times(N-2)$ submatrices defined 
by removing from the $N\times N$ matrix $\big(\Delta + \mu_m {\bf{I}}\big)$ 
the two columns $i$, $j$ and any two rows are invertible may not always be
satisfied. This depends on the topology of the graph and the
choice of the observation nodes $i,j$.
In this subsection, we establish a result for a particular kind of graph 
showing that an arbitrary choice of observation nodes may not satisfy 
the third condition of Theorem \ref{identifiability}. This
result will guide the selection of the observation nodes.
We observed numerically that selecting those nodes following the established 
guidelines leads to fulfill that required condition.

We introduce a joint for a graph.
\begin{definition} A joint is a node whose removal increases
the number of disconnected subgraphs.
\end{definition} 
\begin{figure}[H]
\centerline{
\includegraphics[width=5cm,angle=0]{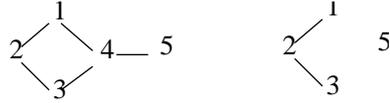}
}
\caption{A connected graph (left) and a graph with two disconnected
subgraphs (right). }
\label{joint5}
\end{figure}
In the example of Fig. \ref{joint5} 4 is a joint of the graph on the
left, since removing it gives a graph with two disconnected subgraphs.
We have the following result.
\begin{theorem} Consider a graph $G$ with a joint at node $k$. Without
loss of generality, we label the nodes $1,2, \dots N$ such that
the sets $\{ 1,2,\dots k-1 \}$ and $\{ k+1,\dots N \}$ correspond to
two disconnected subgraphs of $G$. Choose two indices $i,j \le k$ and
two other indices $p,q \ge k$. Then the  $(N-2)\times(N-2)$ matrix $A_{pq}^{ij}$
obtained by taking out the two columns $i,j$ and the two lines $p,q$ from
$\Delta + \mu_m I$ is not invertible. 
\label{submatrix}
\end{theorem}
{\bf Proof} \\
We assume that $k$ is the joint node.
The matrix $A_m = \Delta + \mu_m I$ has the form
\begin{figure}[H]
\centerline{
\includegraphics[width=8.3cm,angle=0]{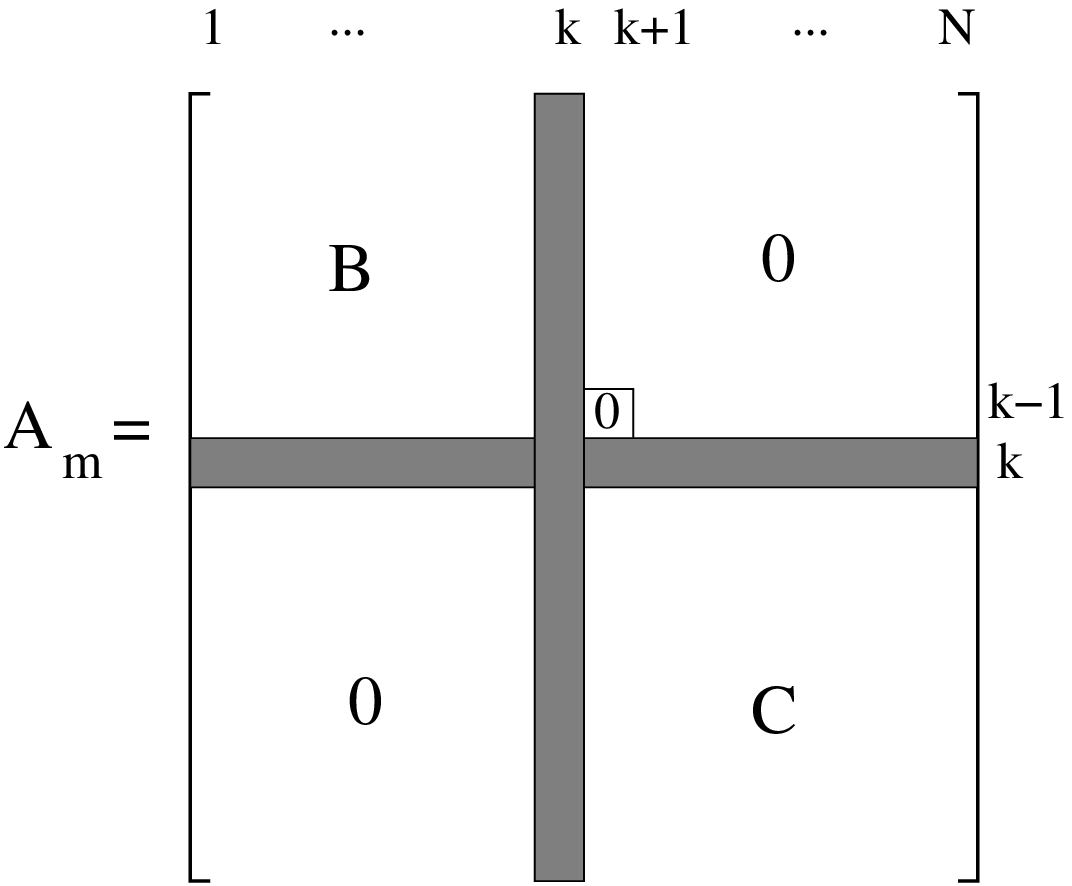}
}
\label{matr}
\end{figure}
There are no links between nodes in $\{ 1,\dots , k-1 \}$
and nodes in $\{ k+1,\dots , N \}$, hence there are two blocks of zeros
in $A_m$.
After removing columns $i$ and $j$ $(i,j \le k)$ and
lines $p$ and $q$ $(p,q \ge k)$ from $A_m$, we get the 
matrix $A_{pq}^{ij}$ 
\begin{figure}[H]
\centerline{
\includegraphics[width=8.3cm,angle=0]{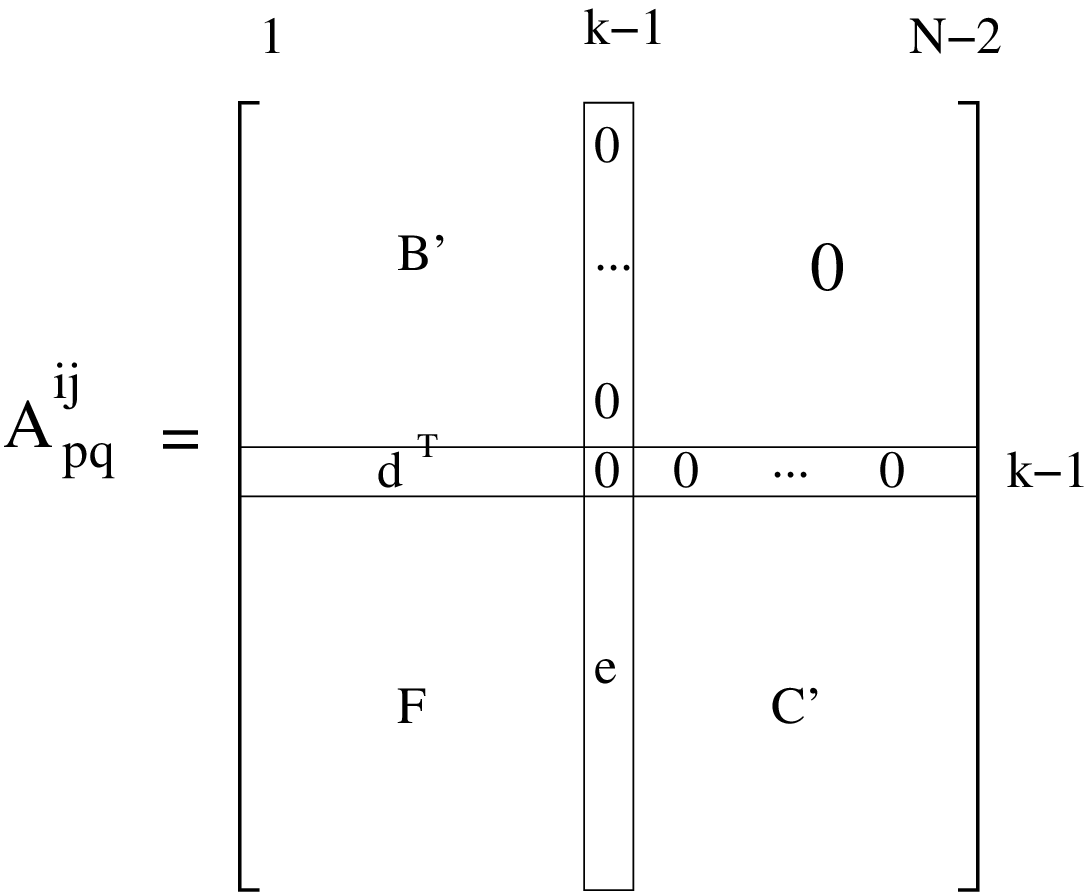}
}
\label{matr2}
\end{figure}
where $d$ is a $1 \times (k-2)$ vector, 
$e$ is a $1 \times (N-k-1)$ vector. The matrix $A_{pq}^{ij}$
is block triangular. The coefficient $0$ on the diagonal of
$A_{pq}^{ij}$ at position $(k-1,k-1)$ comes from 
the coefficient $0$ at position $(k-1,k+1)$ of $A_m$. 
The block triangular matrix $A_{pq}^{ij}$ with a 
$0$ block has a zero determinant and is therefore singular. 
$\blacksquare$

As a consequence of this theorem, the third condition
of the identifiability theorem \ref{identifiability} will
not be satisfied unless we have an observation node
in each part of the graph that will appear as a disconnected subgraph
once a joint is deleted.

Such part of the graph is called a 
maximal bi-connected component in the graph theory literature
\cite{gondran_minoux}. All the joints of a graph can be efficiently
determined using Tarjan's algorithm \cite{tarjan72} in
time $\mathcal{O}(M)$ where $M$ is the number of edges.

\begin{figure}[H]
\centerline{
\includegraphics[width=8.3cm,angle=0]{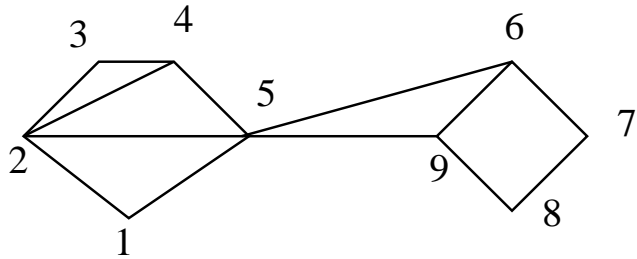}
}
\caption{A graph showing a joint.}
\label{r9a}
\end{figure}

\begin{figure}[H]
\centerline{
\includegraphics[width=5cm,angle=0]{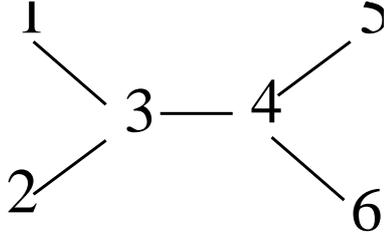}
}
	\caption{A connected graph with two joints.}
\label{joint6}
\end{figure}

An example is the graph shown in Fig. \ref{r9a}, whose node $5$ is a joint.
As a result, putting sensors on any two nodes on the left , say 
$i=1, ~j=2$, will yield a singular matrix in theorem \ref{identifiability}
if we take out two lines in the set $\{ 6, \dots, 9\}$. 
However, placing the two observation nodes such that one node is on 
each side of the joint, for example $i=1$ and $j=7$, fulfills the 
third condition of Theorem \ref{identifiability}.
In Fig. \ref{joint6}, we show another example where nodes 3 and 4 are
joints. Taking out the joints yields four disconnected subgraphs
${1},~ {2},~{5},~{6}$. Then we will need to put sensors on 
nodes $1,2,5,6$. In other
words, there should be a sensor for each maximal bi-connected component
of the graph. 

\section{Identification}

Let $M[\lambda,S]$ be the observation operator introduced in 
(\ref{observations}) defined from recording the state $X$ solution 
of the system (\ref{weq}) in the set of two distinct nodes $\{i,j\}$. 
We consider that the final state $X(T)$ has been determined by the least-square
procedure of section 4.1 so that the vector $P_m$ in (\ref{Notations}) 
is known. Assuming the three conditions of Theorem \ref{identifiability} to 
be satisfied, we focus in this section on developing an identification 
method that determines the two unknown elements $\lambda$ and $S$ of 
the admissible set ${\cal{A}}$ defining the source $\lambda S$ occurring in the system (\ref{weq}). To this end, we proceed as follows: Firstly, we localize the node source position $S$ in the network. Secondly, we propose two different algorithms to identify its emitted signal $\lambda\in L^2(0,T)$.


\subsection{Localization of the node source position $S$}

Since the second assumption of Theorem \ref{identifiability} holds, it follows that all reduced $(N-2)\times(N-2)$ matrix obtained by removing any two rows from the rectangular $N\times(N-2)$ matrix $A_m^{i,j}$ in (\ref{adj1}) is invertible. Therefore, each $(N-2)$ equations of the linear system (\ref{adj1})-(\ref{Pm_ij}) admit a unique solution. Moreover, as the vector $\bar X_m^{i,j}$ defined from the solution $X$ of the system (\ref{weq}) solves the $N$ equations in (\ref{adj1})-(\ref{Pm_ij}) then, all reduced systems of $(N-2)$ equations from (\ref{adj1})-(\ref{Pm_ij}) admit the same unique solution $\bar X^{i,j}_m$.

Hence, to localize the node source position $S$ we proceed as follows: given $m\in \N^*$ such that $\langle \lambda ,\varphi_m\rangle_{L^2(0,T)}\ne 0$, repeat solving a reduced $(N-2)\times(N-2)$ linear system obtained by removing two equations labelled $l_1$ and $l_2$ from the system (\ref{adj1})-(\ref{Pm_ij}) and omitting the source term $\lambda_mS$ i.e., solving only $(N-2)$ equations from the following system:
\begin{eqnarray}
-A_m^{i,j} \bar X^{i,j}_m=P^{i,j}_m
\label{adj2}
\end{eqnarray}

As long as different choices $l_1,l_2$ and $l_1',l_2'$ give the same 
solution $\bar X^{i,j}_m$, this means the source equation i.e., the 
equation containing in its right hand side the unique non-null 
component of the vector $\lambda_m S$, doesn't belong to the 
selected $(N-2)$ equations. However, as soon as a choice $l_1',l_2'$ gives 
a solution different from the one obtained with $l_1,l_2$, then 
the source equation is among the $(N-2)$ solved equations. This process 
leads to determine the label $l$ of the source equation in 
(\ref{adj1}). Moreover, $l$ corresponds to the position of the 
unique non-null component in $S$ which defines the node source.  


\subsection{Identification of the time-dependent signal $\lambda$}

Assuming to be known the position $S$ of the node source $\lambda S$ occurring in the problem (\ref{weq}), we develop two different methods for identifying the unknown time-dependent signal $\lambda\in L^2(0,T)$ emitted by the involved node source: The first method transforms the identification of $\lambda$ into solving a deconvolution problem whereas the second method consists of determining the expansion of $\lambda$ in the orthonormal family of $L^2(0,T)$ made by the normalized eigenfunctions $\{\varphi_m\}$ of the Sturm-Liouville problem introduced in (\ref{Sturm})-(\ref{eigen}).

\subsubsection{\underline{First method}: Deconvolution}

We split the solution $X$ of the system (\ref{weq}) as follows: $X(t)=X^S(t) + X^0(t)$ for all $t\in(0,T)$ where the variables $X^S$ and $X^0$ solve the two following systems:

\begin{eqnarray}
\left\{
\begin{array}{lll}
{\ddot X^S}(t) - \Delta X^S(t)=\lambda(t) S \quad \mbox{in} \; (0,T)\\
X^S(0)=\dot X^S(0)=\vec 0
\end{array}
\right.
\quad \mbox{and} \quad
\left\{
\begin{array}{lll}
{\ddot X^0}(t) - \Delta X^0(t)=0 \quad \mbox{in} \; (0,T)\\
X^0(0)=a \quad \mbox{and} \quad \dot X^0(0)=b
\end{array}
\right.
\label{XSX0}
\end{eqnarray}
{where $a,b$ are exactly the initial conditions of the original problem
\eqref{weq}.}
Using the fundamental solutions obtained in (\ref{Fundamental}) and expanding the variable $X^S$ in the orthonormal family $\{v^1,\dots,v^n\}$ it follows, as in (\ref{X_part}), that 

\begin{eqnarray}
\begin{array}{lll}
X^S(t)=\langle S,v^{1}\rangle \displaystyle \int_0^t \lambda(s)(t-s)ds v^1 + \displaystyle \sum_{n=2}^{N} \frac{\langle S,v^{n}\rangle}{\omega_n}\displaystyle \int_0^t \lambda(s)\sin(\omega_n(t-s))ds v^n,  \quad \mbox{in} \; (0,T)
\end{array}
\label{X_S}
\end{eqnarray}

Afterwards, in view of (\ref{convolution_identifiab}) and according to (\ref{X_S}), it comes that for all $k\in \{1,\dots,N\}$ the $k^{th}$ component $x_k^S$ of the variable $X^S=\big(x_1^S,\dots,x_N^S\big)^{\top}$ is given by

\begin{eqnarray}
\begin{array}{lll}
x_k^S(t)=\displaystyle\int_0^t \lambda(s) \Phi_k(t-s)ds \qquad \forall \; t\in(0,T)\\
\mbox{where:} \quad \Phi_k(t)=\langle S,v^{1}\rangle t v_k^1  + \displaystyle\sum_{n=2}^{N} \frac{\langle S,v^{n}\rangle}{\omega_n} \sin(\omega_nt) v_k^n
\end{array}
\label{X_Sconvolution}
\end{eqnarray}

Therefore, from selecting an observation node $k\in \{i,j\}$ such that $\Phi_k(t)\ne 0$ for almost all $t\in(0,T)$ and since we have: $x_k^S(t)=x_k(t)-x_k^0(t)$ for all $t\in (0,T)$ then, the unknown signal $\lambda\in L^2(0,T)$ occurring in the problem (\ref{weq}) is subject to:

\begin{eqnarray}
\displaystyle\int_0^t\lambda(s)\Phi_k(t-s)ds=x_{k}(t) - x^0_{k}(t) \qquad \forall t\in(0,T)
\label{convolution_ident}
\end{eqnarray}
 
where $x_k^0$ is the $k^{th}$ component of $X^0$. Hence, the identification of $\lambda$ can be transformed into solving the deconvolution problem associated to (\ref{convolution_ident}). Indeed, given a desired number of time steps ${\cal{M}}$, we employ the regularly distributed discrete times: $t_m=m\Delta t$ for $m=0,\dots,{\cal{M}}$ where $\Delta t=T/{\cal{M}}$. Then, using the trapezoidal rule we get

\begin{eqnarray}
\begin{array}{ccccc}
\displaystyle \int_0^{t_{m+1}}\lambda(s)\Phi_k(t_{m+1}-s)ds&=&\displaystyle \sum_{\ell=0}^{m}\int_{t_{\ell}}^{t_{\ell+1}} \lambda(s)\Phi_k(t_{m+1}-s)ds \qquad \qquad\qquad \qquad&\\
&\approx&\displaystyle \frac{\Delta t}{2}\sum_{\ell=0}^{m}\Big(\lambda(t_{\ell}) \Phi_k(t_{m+1-\ell}) + \lambda(t_{\ell+1})\Phi_k(t_{m-\ell})\Big)&\\
&=&\displaystyle \Delta t \sum_{\ell=1}^{m}\lambda(t_{\ell}) \Phi_k(t_{m+1-\ell}) \qquad \qquad\qquad\qquad\qquad &
\end{array}
\label{trapezLam}
\end{eqnarray}

where according to (\ref{X_Sconvolution}), we used $\Phi_k(0)=0$ and assumed that $\lambda(0)=0$. Afterwards, using the notation $\lambda_m\approx \lambda(t_m)$ and in view of (\ref{trapezLam}) we obtain a discrete version of the deconvolution problem associated to (\ref{convolution_ident}) that leads to the following recursive formula:

\begin{eqnarray}
\lambda_m=\displaystyle \frac{1}{\Phi_k(t_1)}\left(\frac{x_{k}(t_{m+1}) - x^0_{k}(t_{m+1})}{\Delta t} - \sum_{\ell=1}^{m-1}\lambda_{\ell} \Phi_k(t_{m+1-\ell})\right), \qquad \forall \; m\ge 1
\label{RecursiveLam}
\end{eqnarray}   

Moreover, by denoting $\Lambda=\big(\lambda_1,\dots,\lambda_{\cal{M}}\big)^{\top}$ and $Q=\big(x_{k}(t_2) - x^0_{k}(t_2),\dots,x_{k}(t_{\cal{M}}) - x^0_{k}(t_{\cal{M}}),x_{k}(t_{\cal{M}}) - x^0_{k}(t_{\cal{M}})\big)^{\top}$ it follows from (\ref{convolution_ident})-(\ref{RecursiveLam}) that the identification of the unknown emitted signal $\lambda$ is transformed into solving the linear system:

\begin{eqnarray}
B\Lambda=\frac{1}{\Delta t} Q
\label{Lin_Syst_lam_n}
\end{eqnarray}

where $B$ is the ${\cal{M}}\times {\cal{M}}$ lower triangular matrix defined for $m=1,\dots,{\cal{M}}$ by

\begin{eqnarray}
\left\{
\begin{array}{lll}
B_{m\ell}=\displaystyle \Phi_k(t_{m+1-\ell}), \qquad\mbox{for} \; \ell=1,\dots,m\\
B_{m\ell}=0, \qquad \mbox{for} \; \ell=m+1,\dots,{\cal{M}}
\end{array} 
\right.
\label{B}
\end{eqnarray}

In practice, the identification of $\Lambda$ using the recursive formula 
in (\ref{RecursiveLam}) which corresponds to the straightforward solution 
of the triangular linear system (\ref{Lin_Syst_lam_n})-(\ref{B}) could 
not yield a stable approximation of the unknown source signal $\lambda$. 
To regularize the system, we employed the Tikhonov 
method that replaces the linear system (\ref{Lin_Syst_lam_n})-(\ref{B}) 
by a penalized least-squares problem under the form:
\begin{eqnarray}
\min_{\Lambda\in \mathbb{R}^{\cal{M}}} \displaystyle \frac{1}{2}\big\|B\Lambda - \frac{1}{\Delta t}Q\big\|_2^2 + \displaystyle \frac{r}{2} \big\|\Lambda\big\|_2^2
\label{Tikhonov} 
\end{eqnarray}
where $r>0$ is the regularization parameter, see for example \cite{Donatelli}. 
A simpler regularization
adding $r$ to the diagonal of the triangular matrix $B$ gave the best results.

\subsubsection{\underline{Second method}: Expansion in a Fourier basis}

This method consists in determining a number of 
coefficients $\lambda_m=\langle \lambda,\varphi_m\rangle_{L^2(0,T)}$ 
defining the expansion of the unknown emitted signal $\lambda$ in 
the orthonormal family $\{\varphi_m\}$ (\ref{Sturm})-(\ref{eigen}). Given 
$m \in \N^*$
a sufficiently large number of $m \in \N^*$ such that the third condition
of Theorem 4.3 is satisfied we can estimate 
\be\label{lambdam} \lambda(t) \approx \sum_m \lambda_m \varphi_m(t) .\ee

The procedure to obtain $\lambda_m$ is the following
\begin{enumerate}
\item Since $S$ is known, select from (\ref{adj1}): $(N-2)$ equations 
that do NOT contain the 
source equation i.e., the equation involving $\lambda_m$ 
in its right-hand side term.  
\item Solve the selected $(N-2)\times(N-2)$ linear system and determine 
$\bar X_m^{i,j}$. 
\item Inject the computed $\bar X_m^{i,j}$ in the 
source equation of (\ref{adj1}) and deduce $\lambda_m$.
\end{enumerate}


\subsection {Algorithm} 

For the clarity of our presentation, we summarize in the 
following algorithm the different steps describing the identification 
method developed in the present paper to determine from the observation 
operator $M[\lambda,S]$ introduced in (\ref{observations}) the two 
unknown elements $\lambda$ and $S$ defining the source $\lambda S$ 
occurring in the problem (\ref{weq}).

{\bf Algorithm}

\begin{enumerate}

\item Determine a {\it{strategic}} set 
of two nodes $\{i,j\}$ satisfying the conditions of 
Theorem \ref{identifiability}. 

\item Compute $X(T)$ by
solving the minimization problem (\ref{Least_squares}) and 
form the rhs of the over-determined linear system (\ref{adj2}).

\item Repeat solving $(N-2)$ equations from 
the linear system (\ref{adj2}) until determining 
the source vector $S$ as explained in subsection 5.1. 

\item Identify the time-dependent 
emitted signal $\lambda$  using one of the two methods suggested
in subsection 5.2. 

\end{enumerate}

\section{Numerical experiments}

To check the effectiveness of the graph source 
identification method developed in the present study, we carry 
out numerical experiments. We first select 
the following $N=5$ nodes graph and choose $1,2$ as our observation
nodes. This will enable us to study the influence of the different 
parameters.
\begin{figure} [H]
\centerline{
\epsfig{file=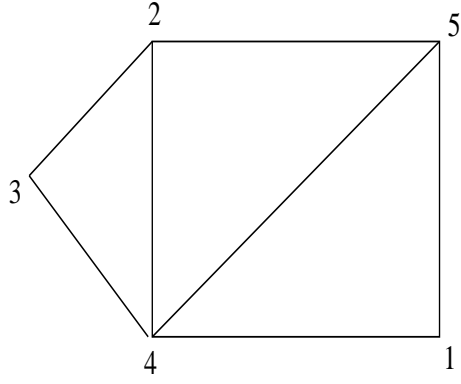,height=5 cm,width=6 cm,angle=0}
}
\caption{5 nodes graph }
\label{g5}
\end{figure}
The associated $5\times 5$ graph Laplacian matrix is given by
\be 
\Delta=
\begin{pmatrix}
-2 & 0 & 0 & 1 & 1 \\
0 & -3 & 1 & 1 & 1 \\
0 & 1 & -2 & 1 & 0 \\
1 & 1 & 1 & -4 & 1 \\
1 & 1 & 0 & 1 & -3
\end{pmatrix}  .
\label{lap_g5}
\ee
Furthermore, the eigenvalues of (\ref{lap_g5}) are: 
$$ 0, ~-3 + \sqrt{2},~-3, ~-3 - \sqrt{2},~-5. $$ 
so that the $\omega_{n},~~n=1,\dots,N$ are
\be
\omega_1 = 0,\quad \omega_2 = 1.259,\quad \omega_3 = 1.732 ,\quad \omega_4 = 2.10 \quad \mbox{and} \quad \omega_5 =2.236
\label{eig_g5}
\ee
The corresponding normalized eigenvectors $v^{n},~n=1,\dots,5$ of (\ref{lap_g5}) are 
\begin{eqnarray}
\label{vec_g5}
v^1 ={1 \over \sqrt{5}}( 1, 1, 1, 1, 1)^T, ~~~
v^2 = {1 \over 2\sqrt{2-\sqrt{2}}}( 1, 1-\sqrt{2}, -1, 0, -1+\sqrt{2})^T, \nonumber \\
v^3 = {1 \over 2} (1, -1, 1, 0, 1)^T,~~~
v^4 = {1 \over 2\sqrt{2+\sqrt{2}}}( 1, 1+\sqrt{2}, -1, 0, -1-\sqrt{2})^T, \\
v^5 = {1 \over 2 \sqrt{5}}( 1, 1, 1, -4, 1)^T \nonumber .
\end{eqnarray}

Besides, according to Definition \ref{strategic} the set of two nodes $\{i=1,j=2\}$ is strategic. In fact, we observe from (\ref{vec_g5}) that $v^{n=1,\dots,N}_i\ne 0$ and $v^{n=1,\dots,N}_j\ne 0$ which implies that the node $i=1$ as well as the node $j=2$ are both strategic. To carry out numerical experiments, we generate synthetic state records by 
solving the forward problem (\ref{weq}) using the time-dependent 
source intensity function defined by
\begin{equation} 
\lambda(t) = \displaystyle\frac{\beta}{2} \left(1+ \tanh(\frac{t-t_L}{w}) - \big(1+ \tanh(\frac{t-t_R}{w})\big) \right) ,
\label{lambda_tanh}
\end{equation}
for all $t\in(0,T_0)$ whereas $\lambda(t)=0$ for all $t\ge T_0$. 
We set the observation operator $M[\lambda,S]$ at the strategic set of 
two nodes $\{i=1,j=2\}$ as in (\ref{observations}). To this end, we used in 
(\ref{lambda_tanh}) the coefficients: 
\be\label{coeff_g5} 
\beta=3,~~t_L=0.3 T, w=0.01 T,~~ t_R=0.6T ,~~ T=100,~~s=3. \ee
Then the time active limit $T_0$ in (\ref{Ad}) is $T_0=70$. Moreover, given 
a desired number of time steps ${\cal{M}}$, we employ the regularly 
distributed discrete times $t_m=m\Delta t$ for $m=0,\dots,{\cal{M}}$ where 
$\Delta t=T/{\cal{M}}$ and solve the forward problem 
(\ref{weq}), (\ref{lambda_tanh}) using a variable step Runge-Kutta solver of
order 4-5.

First, we consider the first step of the {\bf{Algorithm.}}, the estimation
of the final state $X(T)$ using least squares.

\subsection{ Least squares fit }

The two nodes $i=1$ and $j=2$ defining the observation operator 
$M[\lambda,S]$ are both strategic. Then it follows from Remark 
\ref{remark_ls} that using only the state records $x_k(t)$ for $t\in(T^*,T)$ 
and $k=1$ or $k=2$ yields uniqueness of the sought final state vector 
$X(T)$. Then, using the regularly distributed discrete times: 
$t_m=m\Delta t$ for $m=0,\dots,{\cal{M}}$ where $\Delta t=T/{\cal{M}}$ 
and assuming there exists ${\cal{M}}^*$ such that 
$T^*={\cal{M}}^* \Delta t$, we sample the observed data 
$d_k(t)\equiv x_k(t)$ 
as $d_k(t_m)$ for $m=1,\dots,{\cal{M}}$. Therefore, we aim to identify 
the final state vector $X(T) =\big(x_1(T),x_2(T), \dots, x_N(T)\big)^{\top}$ 
by solving the least squares problem (\ref{Least_squares}).

We generate synthetic state records $D_k$ by solving the problem 
(\ref{weq}) for a source located at node $s=3$ 
wth the time-dependent intensity function $\lambda$ introduced in 
(\ref{lambda_tanh}). For simplicity, we set in the least squares problem 
(\ref{Least_squares}): ${\cal{M}}=100$ and $T^*=T^0$ 
which implies ${\cal{M}}^*=70$. The relative error $r_k$ between 
the generated (exact) final state vector $X_{ex}(T)$ and the identified 
final state vector $X^k_{id}(T)$ 
\be
r_k =\displaystyle\frac{\big\|X^k_{id}(T)- X_{ex}(T)\big\||_2}{\big\|X_{ex}(T)\big\|_2} ,
\label{rel_err}
\ee
is presented in Table \ref{tab1} for different observation nodes $k$
in single and double precision.

\begin{table} [H]
\centering
\begin{tabular}{|l|c|c|c|r|}
   \hline\hline
Observation node     &    $k=1$    &    $k=2$   &   $k=3$    &  $k=5$ \\ \hline
$r_k$ double& $9.87~10^{-06}$ & $1.08~10^{-05}$  & $ 1.46~ 10^{-05}$ & $9.74~10^{-06}$ \\ \hline
$r_k$ single& $4~10^{-05}$  &$3.59~10^{-05}$ & $4.19~10^{-05}$  & $4~10^{-05}$  \\ 
\hline\hline
\end{tabular}
\caption{Relative error $r_k$ in the least-square estimation of $X(T)$
for different observed nodes $k=1,2,3,5$ in single and double precision.}
\label{tab1}
\end{table}
The results presented in Table \ref{tab1} confirm our claim 
in Remark \ref{remark_ls}; they show that the state records at 
a single strategic node $k$ lead to identify with accuracy the 
sought final state vector $X(T)$.

However, if the node $k$ is not strategic, we show with the following
example that the solution of the non regularized least-square problem 
from (\ref{Least_squares}) may not exist or not be unique.
Consider the node $k=4$ which is not strategic
because from (\ref{vec_g5}) $v^2_4=v^3_4=v^4_4=0$. 
The least square consists in minimizing with respect to the parameters 
$y_{n}(T), ~~\dot y_{n}(T),~~n=1,\dots,N$ 
the following function:
\be 
E_k\big(y_{n=1,\dots,N}(T),\dot y_{n=1,\dots,N}(T)\big)=
\displaystyle\sum_{m={\cal{M}}^*}^{{\cal{M}}} ( x_k(t_m)- d_k(t_m))^2,
\label{error}
\ee
where $x_k$ is defined as 
in (\ref{Express_X_ident}) by
\be 
x_k(t)=\Big(y_1(T) +(t-T)\dot y_1(T)\Big)v_k^1
+ \sum_{n=2}^N \Big(y_n(T) \cos (\omega_n (t-T)) + 
{\dot y_n(T) \over \omega_n} \sin (\omega_n (t-T))\Big)v_k^n ,
\label{x1}
\ee

The function $E_k$ in (\ref{error}) can be written as
\be 
E_k(Y)=\big\|A_k Y-D_k\|_2^2 ,
\label{error_matrix}
\ee
where $Y=\big(y_1(T),\dot y_1(T), \dots,y_N(T),\dot y_N(T)\big)^{\top}$, $D_k=\big(d_k(t_{{\cal{M}}^*}), \dots d_k(t_{\cal{M}})\big)^{\top}$ and
{\tiny 
\begin{eqnarray}
A_k= \begin{pmatrix}
v_k^1 & (t_{{\cal{M}}^*}-T)v_k^1 & \cos (\omega_2 (t_{{\cal{M}}^*}-T))v_k^2 & { \sin (\omega_2 (t_{{\cal{M}}^*}-T)) \over \omega_2}v_k^2& \dots & \cos (\omega_N (t_{{\cal{M}}^*}-T))v_k^N & {\sin (\omega_N (t_{{\cal{M}}^*}-T)) \over \omega_N}v_k^N\\
v_k^1  & (t_{{\cal{M}}^* +1}-T)v_k^1  & \cos (\omega_2 (t_{{\cal{M}}^* +1}-T))v_k^2 & { \sin (\omega_2 (t_{{\cal{M}}^* +1}-T))\over \omega_2}v_k^2 &\dots & \cos (\omega_N (t_{{\cal{M}}^* +1}-T))v_k^N & {\sin (\omega_N (t_{{\cal{M}}^* +1}-T)) \over \omega_N}v_k^N \\
\vdots & \vdots & \vdots & \vdots &  \vdots & \vdots &  \\
v_k^1 & (t_{\cal{M}}-T)v_k^1 & \cos (\omega_2 (t_{\cal{M}}-T))v_k^2 & { \sin (\omega_2 (t_{\cal{M}}-T)) \over \omega_2}v_k^2 & \dots & \cos (\omega_N (t_{\cal{M}}-T))v_k^N &  {\sin (\omega_N (t_{\cal{M}}-T)) \over \omega_N} v_k^N
\end{pmatrix} .
\label{matrix_Ak}
\end{eqnarray}
}
For $k=4$, the matrix $A_{k=4}$ is exactly singular because
it has three columns of zeros, corresponding 
to $v_4^2, v_4^3 $ and $v_4^4$. Then the linear
system has no solution or an infinite number of solutions
and no information can be obtained 
on $y_2, \dot y_2,y_3, \dot y_3,y_4, \dot y_4$.
This means that the least square problem for the $k=4$ non strategic
node cannot be solved.

\subsection{Localization of the unknown node source} 

We now assume the final state vector $X(T)$  to be known and focus 
on the localization of the unknown node source using the time records 
in $(0,T)$ of the state $x_1(t)$ and $x_2(t)$ at the two 
strategic nodes $\{i=1,j=2\}$.

The adjoint problem (\ref{adj1}) 
is the following linear system of $N=5$ equations and $N-2=3$ unknowns:
\be
-A^{1,2}  \bar X_m^{1,2}= P_m^{1,2} + \lambda_m S \label{id_sys}
\ee
where $\bar X_m^{1,2} \in \R^{N-2}$ and the coefficient $P_m^{1,2}$ is
given by (\ref{Pm_ij}).

We carried out numerical experiments on the localization of the node 
source from state records $x_1(t)$ and $x_2(t)$ for $t\in(0,T)$ 
generated by a source occurring at node $s=3$ i.e., 
$S=\big(0,0,1,0,0\big)^{\top}$ forcing the graph with the time-dependent 
intensity function $\lambda$ introduced in (\ref{lambda_tanh}). We applied 
the procedure described above i.e., repeat taking out two equations 
$l_1$ and $l_2$ from the $N$ equations in (\ref{id_sys})
and solve the remaining $(N-2)$ equations. The results are presented 
in the following table where for each couple $(l_1, l_2)$ 
we give the solution $\bar X_m^{1,2}$ together with the
error between the computed solution $\bar X_m^{1,2}$ and 
the exact solution $\bar X_{ex}^{1,2}$ generated by solving 
the direct problem (\ref{weq}) with $\lambda$ and $S$.
\begin{table} [H]
\centering
\begin{tabular}{|l|c|c|r|}
   \hline
 &  &    &   \\
$l_1$  & $l_2$  & $\bar X_m^{1,2}$  & $\|\bar X_m^{1,2} -\bar X^{1,2}_{ex}\|_2$ \\ 
 &  &    &   \\ \hline
1& 2&    $(2020.704778~ , 2020.72667~ ,2020.732166~)^T$ & $3.58 ~ 10^{-3}$ \\ \hline
1& 4&    $(2023.289145~ , 2025.89286~ ,2022.454794~)^T$  & $2.75 ~ 10^{-3}$ \\ \hline
1& 5&    $(2021.652750~ ,2022.621687~, 2027.362363~)^T$ & $4.19 ~ 10^{-3}$ \\ \hline
2& 4&    $(2021.652808~, 2022.621803~, 2021.364083~)^T$ & $2.97 ~ 10^{-3}$ \\ \hline
2& 5&    $(2020.985776~, 2021.288396~, 2022.697489~)^T$ & $3.43 ~ 10^{-3}$ \\ \hline
4& 5&    $(2027.650916~, 2034.612098~, 2009.373787~)^T$ & $8.39 ~ 10^{-3}$ \\ \hline
 &  &    &   \\  \hline
1& 3&    $(2027.651070~, 2022.621687~, 2021.364044~)^T$ & $1.64 ~ 10^{-7}$  \\ \hline
2& 3&    $(2027.651492~, 2022.621803~, 2021.364083~)^T$ & $3.72 ~ 10^{-7}$  \\ \hline
3& 4&    $(2027.650916~, 2022.621803~, 2021.364083~)^T$ & $1.28 ~ 10^{-7}$  \\ \hline
3& 5&    $(2027.650916~, 2022.621687~, 2021.364198~)^T$ & $8.52 ~ 10^{-8}$ \\
\hline
\end{tabular}
\caption{Identification of a source at node $s=3$:  solution
$(\bar X_3, \bar X_4  ,\bar X_5)$ of the over-determined system
(\ref{adj2}) when lines $l_1,l_2$ are taken out. 
The calculations were done with 64 bit arithmetics and the
parameters is $m=1$.
}
\label{tab2}
\end{table}

The numerical results presented in Table \ref{tab2} show that 
the source localization procedure developed in the present study 
leads to identify with accuracy the node source generating 
the state records $x_1(t),x_2(t)$ for $t\in(0,T)$. Indeed, 
the first part in Table \ref{tab2} shows that as long as 
the source equation $s=3$ is not among the two equations 
$l_1$ and $l_2$ taken out from the $N$ equations in (\ref{id_sys}), the 
solutions of the remaining $(N-2)$ equations are different and the 
error, about $10^{-3}$ is relatively important.
However, when the source equation 
$s=3$ is among the two lines $l_1$ or $l_2$ as in the second part of 
Table \ref{tab2}, we observe that all remaining $(N-2)$ equations 
admit about the same solution $\bar X^{1,2}_m$ that 
fits closely the exact solution $\bar X^{1,2}_{ex}$ because the error 
is about $10^{-5}$. Similar calculations done with 32 bits arithmetic
yield an error of $10^{-5}$ when $l_1=s$ or $l_2=s$, otherwise
the error is about $10^{-3}$, see Table \ref{tab3}.

The results were obtained for a source at node $s=3$. Placing
the source at the other nodes of the $5$-nodes graph yielded similar
results to those presented in Table \ref{tab2}.

We now examine the influence of the order $m$ of the 
eigenfunction $\varphi_m$ needed to form the adjoint problem.
We repeated the calculations shown in Table \ref{tab2} for
$m=5$ and $m=10$ and the results are shown in tables 
\ref{tab3} and \ref{tab4} respectively.
\begin{table} [H]
\centering
\begin{tabular}{|l|c|c|c|c|r|}
   \hline\hline
       &       &                  &  \\ 
$l_1$  & $l_2$  & $\bar X_m^{1,2}$  & $\|\bar X_m^{1,2} -\bar X^{1,2}_{ex}\|_2$ \\ 
       &       &                  &  \\ \hline
 1     & 2     &$\big(871.11,~871.12,~871.12\big)^{\top}$ & $1.54 ~ 10^{-3}$\\\hline
 1     & 4     &$\big(871.62,~872.11,~871.45\big)^{\top}$ & $1.17 ~ 10^{-3}$\\\hline
 1     & 5     &$\big(871.3,~871.49,~872.4\big)^{\top}$   & $1.82 ~ 10^{-3}$\\\hline
 2     & 4     &$\big(871.32,~871.53,~871.26\big)^{\top}$ & $1.24 ~ 10^{-3}$\\\hline
 2     & 5     &$\big(871.18,~871.24,~871.55\big)^{\top}$ & $1.47 ~ 10^{-3}$\\\hline
 4     & 5     &$\big(872.4,~873.66,~869.13\big)^{\top}$  & $3.49 ~ 10^{-3}$\\\hline
       &       &                                       &  \\ \hline
 1     & 3     &$\big(872.46,~871.49,~871.24\big)^{\top}$ & $7.04 ~ 10^{-5}$\\\hline
 2     & 3    &$\big(872.62,~871.53,~871.26\big)^{\top}$ & $2.57 ~ 10^{-4}$\\\hline 
 3     & 4     &$\big(872.4,~871.53,~871.26\big)^{\top}$ & $5.10 ~ 10^{-5}$\\\hline
 3     & 5     &$\big(872.4,~871.49,~871.3\big)^{\top}$  & $5.01 ~ 10^{-5}$\\
\hline\hline
\end{tabular}
\caption{Localization of the node source $s=3$ with $m=5$; the other
parameters are the same as in Table \ref{tab2}.}
\label{tab3}
\end{table}

\begin{table} [H]
\centering
\begin{tabular}{|l|c|c|c|c|r|}
   \hline\hline
       &       &                  &  \\
$l_1$  & $l_2$  & $\bar X_m^{1,2}$  & $\|\bar X_m^{1,2} -\bar X^{1,2}_{ex}\|_2$ \\ 
       &       &                  &  \\ \hline
1      & 2     &$\big(-439.08,~-439.24,~-439.24\big)^{\top}$&$2.49~ 10^{-3}$\\\hline
1      & 4     &$\big(-439.67,~-440.36,~-439.63\big)^{\top}$&$1.83~ 10^{-3}$\\\hline
1      & 5     &$\big(-439.31,~-439.67,~-440.68\big)^{\top}$&$3.19~ 10^{-3}$\\\hline
2      & 4     &$\big(-439.48,~-440,~-439.5\big)^{\top}$&$1.51~ 10^{-3}$\\\hline
2      & 5    &$\big(-439.21,~-439.48,~-440.03\big)^{\top}$ &$2.33~ 10^{-3}$\\\hline
4      & 5   &$\big(-440.16,~-441.28,~-438.22\big)^{\top}$ & $4.67~ 10^{-3}$\\\hline
       &      &                                        &  \\ \hline
1      & 3   &$\big(-440.6,~-439.67,~-439.39\big)^{\top}$& $1.26~ 10^{-3}$\\\hline
2      & 3   &$\big(-441.78,~-440,~-439.5\big)^{\top}$ &   $3.98~ 10^{-3}$\\\hline
3      & 4   &$\big(-440.16,~-440,~-439.5\big)^{\top}$ &  $7.94~ 10^{-4}$ \\ \hline 
3     &  5   &$\big(-440.16,~-439.67,~-439.83\big)^{\top}$ &   $8.04~ 10^{-4}$ \\ \hline 
\hline\hline
\end{tabular}
\caption{Localization of the node source $s=3$ with $m=10$; the other
parameters are the same as in Table \ref{tab2}.}
\label{tab4}
\end{table}
Compared to Table \ref{tab2} the 
numerical results of Table \ref{tab3} show that the error
between the exact solution and the identified $\bar X_m^{12}$
has increased from $10^{-7}$ to $10^{-5}$.  For $m=10$
shown in Table \ref{tab4}, there is no significant difference
between the error when the source line $s=3$ is included and
when it is not.
This suggests to use small values of $m$ for the 
source localization procedure.
In fact, a large $m$ corresponds to a higher 
frequency eigenfunction $\varphi_m$ from (\ref{eigen}) so that
the integrands defining the terms
involved in the right hand side $P_m^{i,j}$ in
(\ref{adj2})-(\ref{Pm_ij}) are highly oscillatory. The integrals
become therefore less accurate.

To conclude this section, note that to identify the unknown source node, 
we developed a simple algorithm {where 
for all lines $l_1$, we pick two other lines $l_2$ and $l_3$.
For the two couples of lines, $(l_1,l_2)$ and $(l_1,l_3)$, 
we extract the two lines of \eqref{adj1} and
solve the reduced system. If the norm of the difference of
the solutions is smaller than a threshold, we retain $l_1$
as the node source. This yields a complexity $\mathcal{O}(N^4)$ ,
$N$ for the sweep in $l_1$ and $N^3$ for solving the reduced  
linear system.}
We used this naive algorithm successfully in preliminary blind 
tests -where the source was not known to the user- on 
a $14$-node network.

\subsection{Identification of the time-dependent source signal $\lambda$}

In this last part of our numerical experiments, we assume the source position $S$ to be now known and aim to identify its time-dependent intensity function $\lambda$ using the two developed methods: Deconvolution and expansion in the complete orthonormal family $\{\varphi_m\}_m$. For both methods, we use the same state records $x_1(t)$ and $x_2(t)$ for $t\in(0,T)$ generated by a source occurring 
in node $s=3$ i.e., $S=\big(0,0,1,0,0\big)^{\top}$ forcing the graph with the time-dependent intensity function $\lambda$ introduced in (\ref{lambda_tanh}). In addition, we use the already introduced regularly distributed discrete times: $t_m=m\Delta t$ for $m=0,\dots,{\cal{M}}$ where $\Delta t=T/{\cal{M}}$ for ${\cal{M}}=100$.

\subsubsection{Identification of $\lambda$ using deconvolution}

This first method identifies the time-dependent source intensity 
function $\lambda$ by solving the discrete version of the deconvolution 
problem obtained in (\ref{convolution_ident}). This leads to 
determine $\lambda^m\approx \lambda(t_m)$ for $m=1,\dots,{\cal{M}}$ 
from the recursive formula obtained in (\ref{RecursiveLam}). This
formula is equivalent to solve the
following ${\cal{M}}\times {\cal{M}}$ triangular linear system 
\be
\begin{pmatrix}
\Phi_k(t_1)   &      0         & \dots      &     &         &0 \\
\Phi_k(t_2)   &   \Phi_k(t_1)   & 0  &\dots    &         & 0\\
\Phi_k(t_3)   &   \Phi_k(t_2)   & \Phi_k(t_1) & 0     & \dots        & 0\\
\vdots      &    \vdots     & &    &   \vdots  & \\
\Phi_k(t_{\cal{M}})   & \Phi_k(t_{{\cal{M}}-1}) & \Phi_k(t_{{\cal{M}}-2}) & \Phi_k(t_{{\cal{M}}-3}) & \dots  & \Phi_k(t_1) 
\end{pmatrix}
\begin{pmatrix}
\lambda^1 \\
\lambda^2 \\
\lambda^3\\
\vdots \\
\lambda^{\cal{M}}
\end{pmatrix}
=
\begin{pmatrix}
b_k^2\\
b_k^3\\
b_k^4\\
\vdots \\
b_k^{{\cal{M}}+1}
\end{pmatrix}
\label{LS_discrete}
\ee
where the right hand side is given by the state record 
$x_k(t)$ for $t\in(0,T)$ where $k=1$ or $k=2$, 
$b_k^m=\big(x_k(t_m) - x_k^0(t_m)\big)/{\Delta t}$ 
(from (\ref{convolution_ident})) and where
$\Phi_k(t_m)$ is given by (\ref{X_Sconvolution}).

We carry out numerical experiments on the identification of the 
time-dependent signal $\lambda$ emitted by the already localized 
node source $S$ using the deconvolution method developed in 
(\ref{RecursiveLam})-(\ref{Tikhonov}). For each experiment, we 
generate simulated measures by solving the system (\ref{weq}) with 
a different kind of signal $\lambda=\lambda_l, ~~l\in\{1,2,3\}$ 
where $\lambda_1$ 
is given by (\ref{lambda_tanh}) and (\ref{coeff_g5}) and 
\begin{eqnarray}
\begin{array}{lll}
\bullet \; \lambda_2(t)=\sin(\pi t/T^0), \qquad \quad\;\; \forall t\in(0,T^0) \qquad \mbox{and} \qquad \lambda_2(t)=0,\quad \forall t\ge T^0\\
\bullet \; \lambda_3(t)=\displaystyle \sum_{n=1}^3 c_ne^{-\alpha_n(t-\tau_n)^2}, \quad \forall t\in(0,T^0) \qquad \mbox{and} \qquad \lambda_3(t)=0,\quad \forall t\ge T^0
\end{array}
\label{lambdas}
\end{eqnarray} 
where $\beta=2$ and $c_1=1.2$, $c_2=0.4$, $c_3=0.6$, $\alpha_1=10^{-6}$, 
$\alpha_2=5\times 10^{-5}$, $\alpha_3=10^{-6}$ and 
$\tau_1=4500$, $\tau_2=6500$, $\tau_3=8500$. 
We use here a large time interval and the same sampling rate as in
the previous sections, 
$$T=14400, ~~T^0=10800~~~ {\rm and} ~~~{\cal{M}}=14400  .$$ 
{Note that here the typical periods of oscillation of the network
as given by (\ref{eig_g5}) are much smaller than $T$.}

\begin{figure}[H]
\begin{center}
\begin{eqnarray*}
\begin{array}{cccc}
{\bf{(a)}}&&{\bf{(b)}}\\
\includegraphics[width=7cm,angle=0]{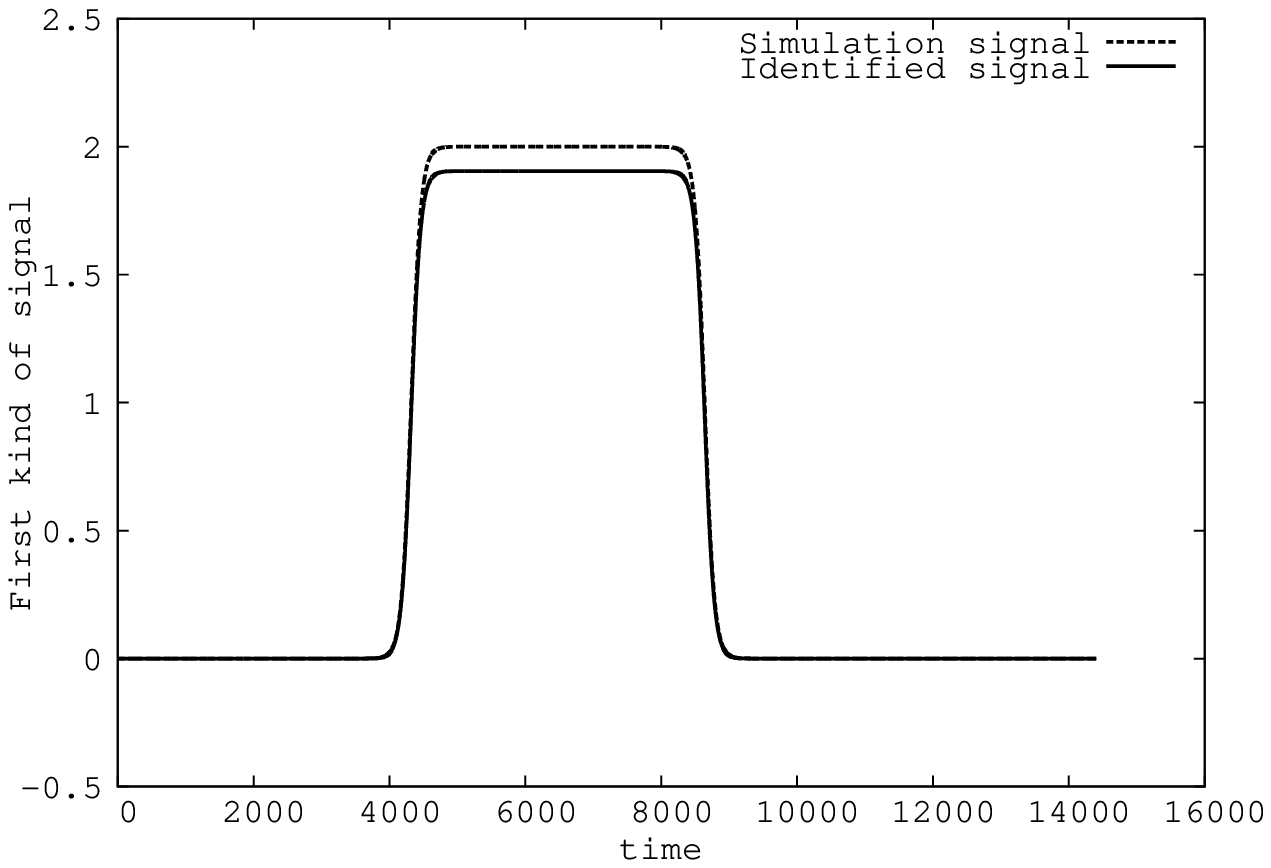}&&\includegraphics[width=7cm,angle=0]{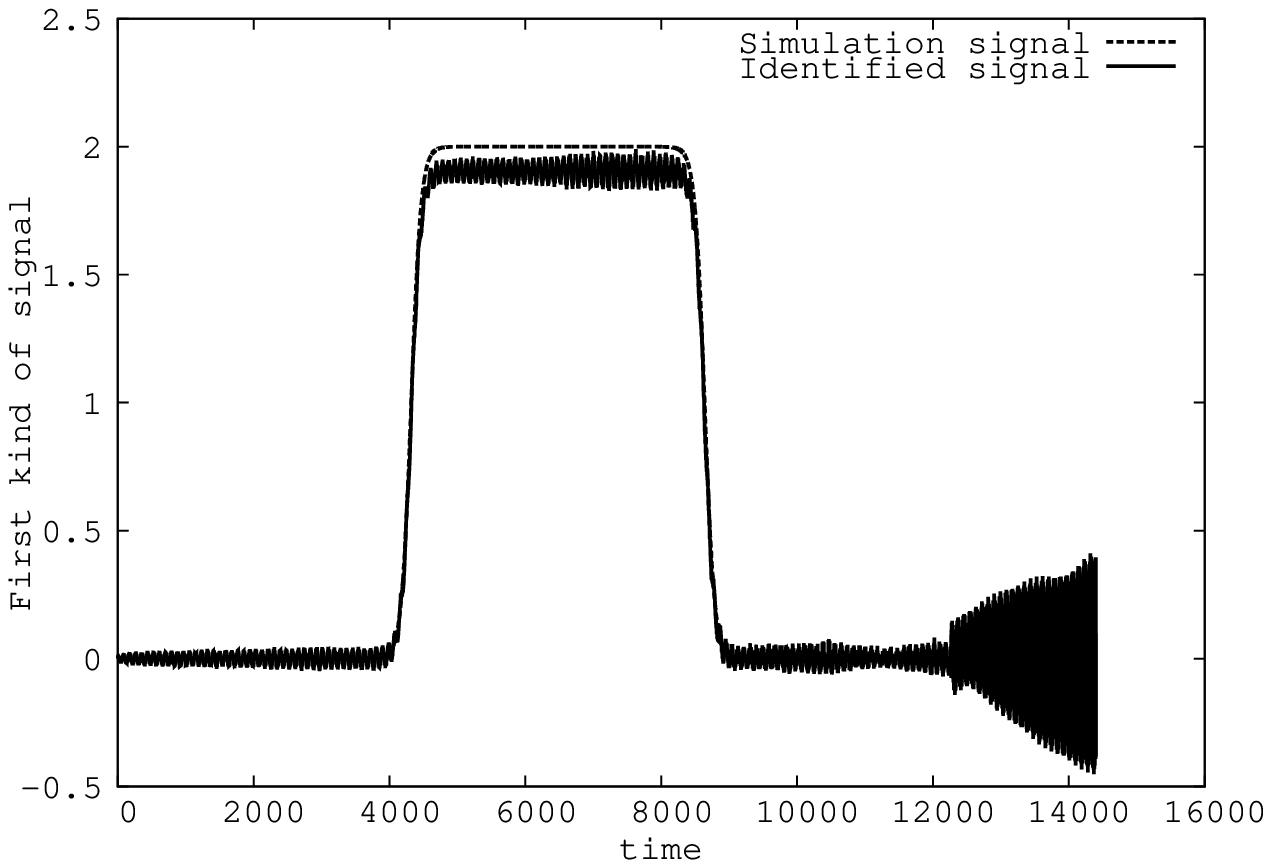}
\end{array}
\end{eqnarray*}
\caption{{\bf{(a)}} Gaussian Noise $0\%$, $\mbox{Error}_1=0,04$ \hspace{0.4cm} {\bf{(b)}} Gaussian Noise $1\%$, $\mbox{Error}_1=0,08$}
\label{Fig1_lam1}
\end{center}
\end{figure}

\begin{figure}[H]
\begin{center}
\begin{eqnarray*}
\begin{array}{cccc}
{\bf{(c)}}&&{\bf{(d)}}\\
\includegraphics[width=7cm,angle=0]{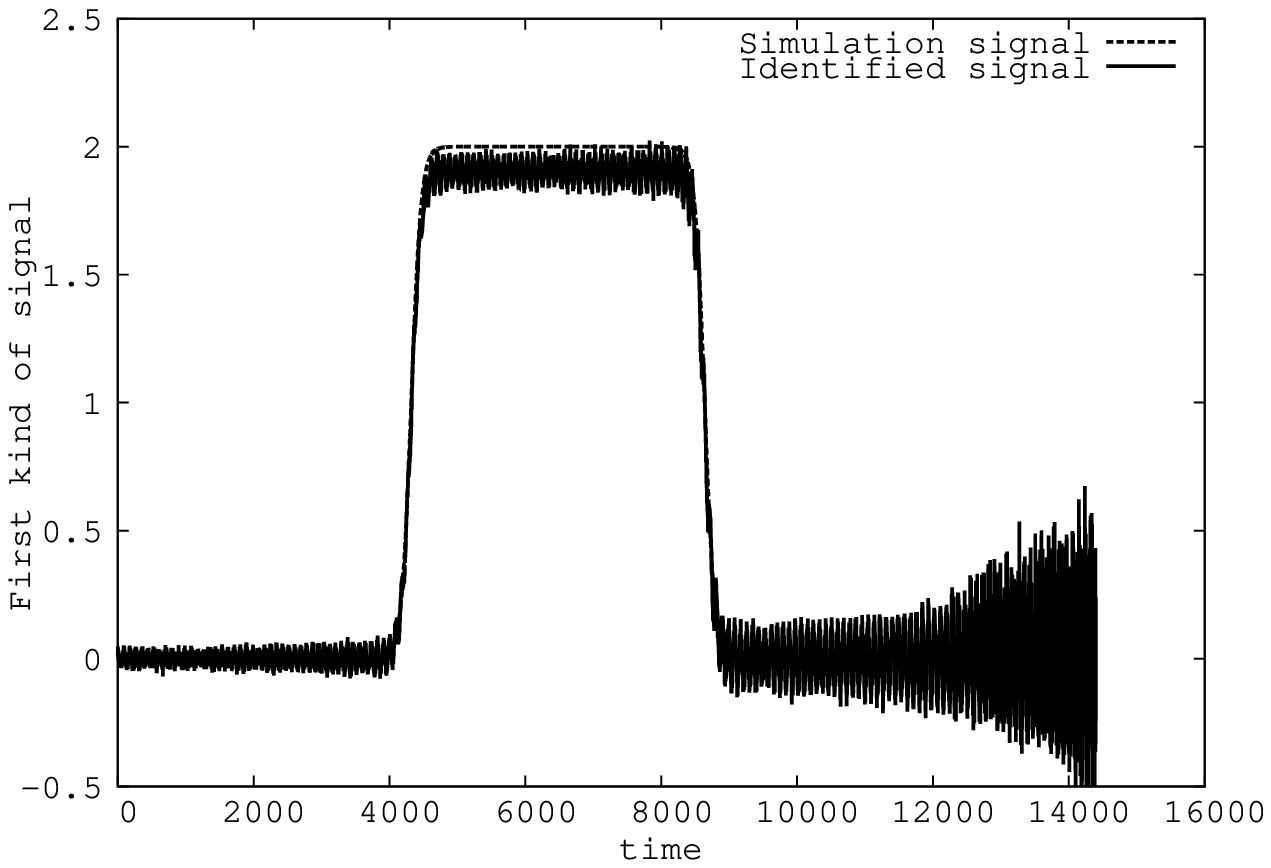}&&\includegraphics[width=7cm,angle=0]{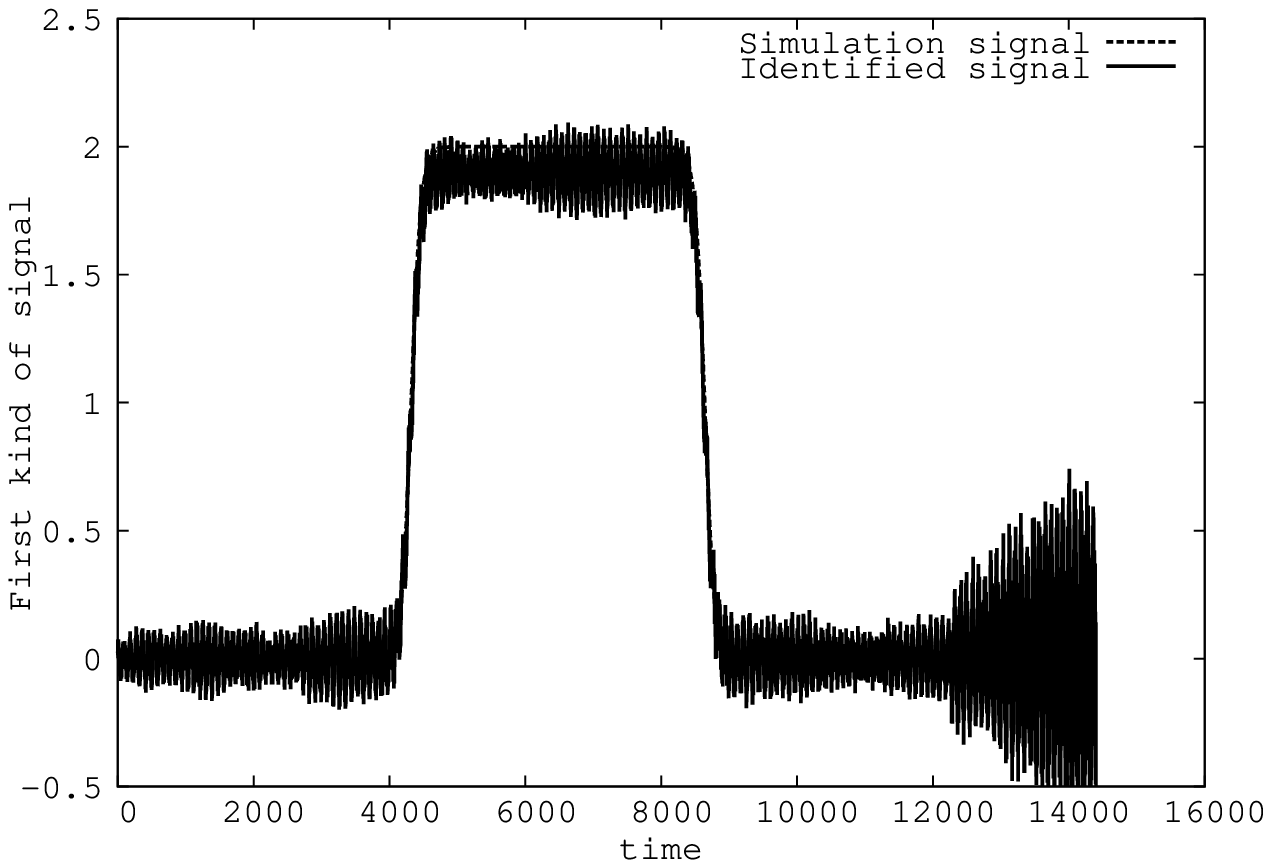}
\end{array}
\end{eqnarray*}
\caption{{\bf{(c)}} Gaussian Noise $3\%$, $\mbox{Error}_1=0,09$ \hspace{0.4cm} {\bf{(d)}} Gaussian Noise $5\%$, $\mbox{Error}_1=0,11$}
\label{Fig2_lam1}
\end{center}
\end{figure}

We start by presenting the numerical results on the identification of the first kind of signal in (\ref{lambdas}) i.e., applying the deconvolution method using simulated measures generated by $\lambda=\lambda_1$ and emitted by a source located in the node $3$ i.e., $S=(0,0,1,0,0)^{\top}$. First, we present 
in Figures \ref{Fig1_lam1} and \ref{Fig2_lam1} the curve, labeled 
"Simulation signal", of $\lambda_1$ introduced in (\ref{lambda_tanh}) 
as well as the curve of the identified signal using measures taken 
at the observation node $k=1$ and subject to different intensities of a 
Gaussian noise. 
Second, we identify $\lambda_1$ using measures 
taken for each experiment in a different observation node 
$k\in\{2,3,4,5\}$ and subject to a Gaussian noise of intensity $3\%$. The 
results are presented in Table \ref{tab_lam1} where we indicate 
for each observation node $k$ the corresponding $L^2$ relative error:
\begin{eqnarray}
\mbox{Error}_l=\displaystyle\frac{\|\Lambda_l-\Lambda_l^{ident}\|_2}{\|\Lambda_l\|_2}
\label{error}
\end{eqnarray}
where $\|\cdot\|_2$ designates the euclidean norm, $\Lambda_l=\big(\lambda_l(t_1),\dots,\lambda_l(t_{\cal{M}})\big)$ and $\Lambda_l^{ident}=\big(\lambda^1_l,\dots,\lambda_l^{\cal{M}}\big)$ with $\lambda_l^i$ is the identified value of $\lambda_l(t_i)$.


The numerical results presented in Figures \ref{Fig1_lam1} 
and \ref{Fig2_lam1} show that the identification method based on 
deconvolution developed in the present paper leads to identify the 
emitted signal $\lambda_1$ and is relatively stable with respect to 
the introduction of a noise on the measures. Then, we give in the 
following table the $L^2$ relative error computed from (\ref{error}) 
between the "Simulation signal" and the identified signal while 
using measures taken at a different observation node $k$:  
\begin{table} [H]
\centering
\begin{tabular}{|l|c|c|c|c|c|r|}
   \hline\hline
Observor node     &    $k=1$    &    $k=2$   &   $k=3$    &  $k=4$ & $k=5$ \\ \hline
$\mbox{Error}_1$ & $0,093$ & $0,086$  & $0,072$ & - & $0,096$ \\  
\hline\hline
\end{tabular}
\caption{Identification of $\lambda_1$: Source $S=(0,0,1,0,0)^{\top}$ and Gaussian Noise $3\%$}
\label{tab_lam1}
\end{table}

\begin{figure} [H]
\begin{center}
\begin{eqnarray*}
\begin{array}{cccc}
{\bf{(a)}}&&{\bf{(b)}}\\
\includegraphics[width=7cm,angle=0]{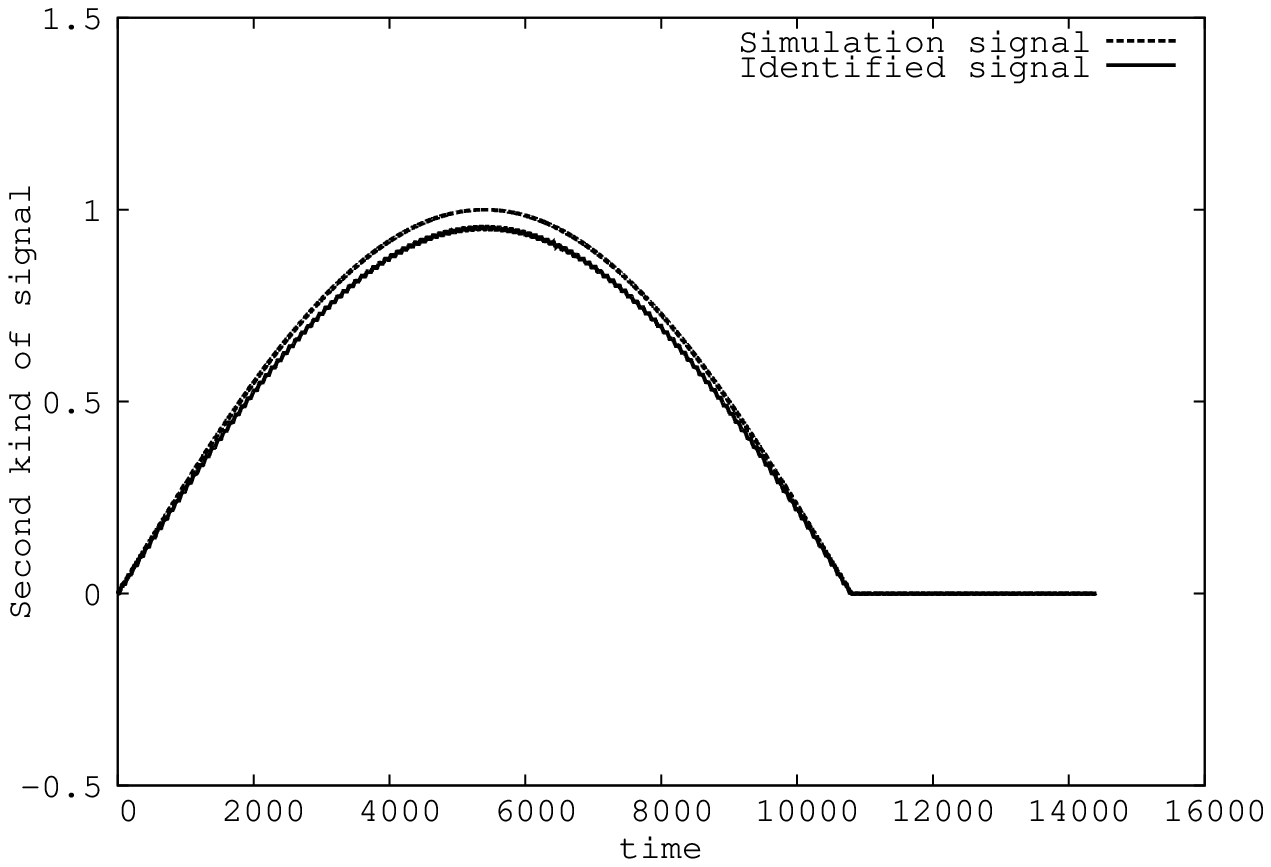}&&\includegraphics[width=7cm,angle=0]{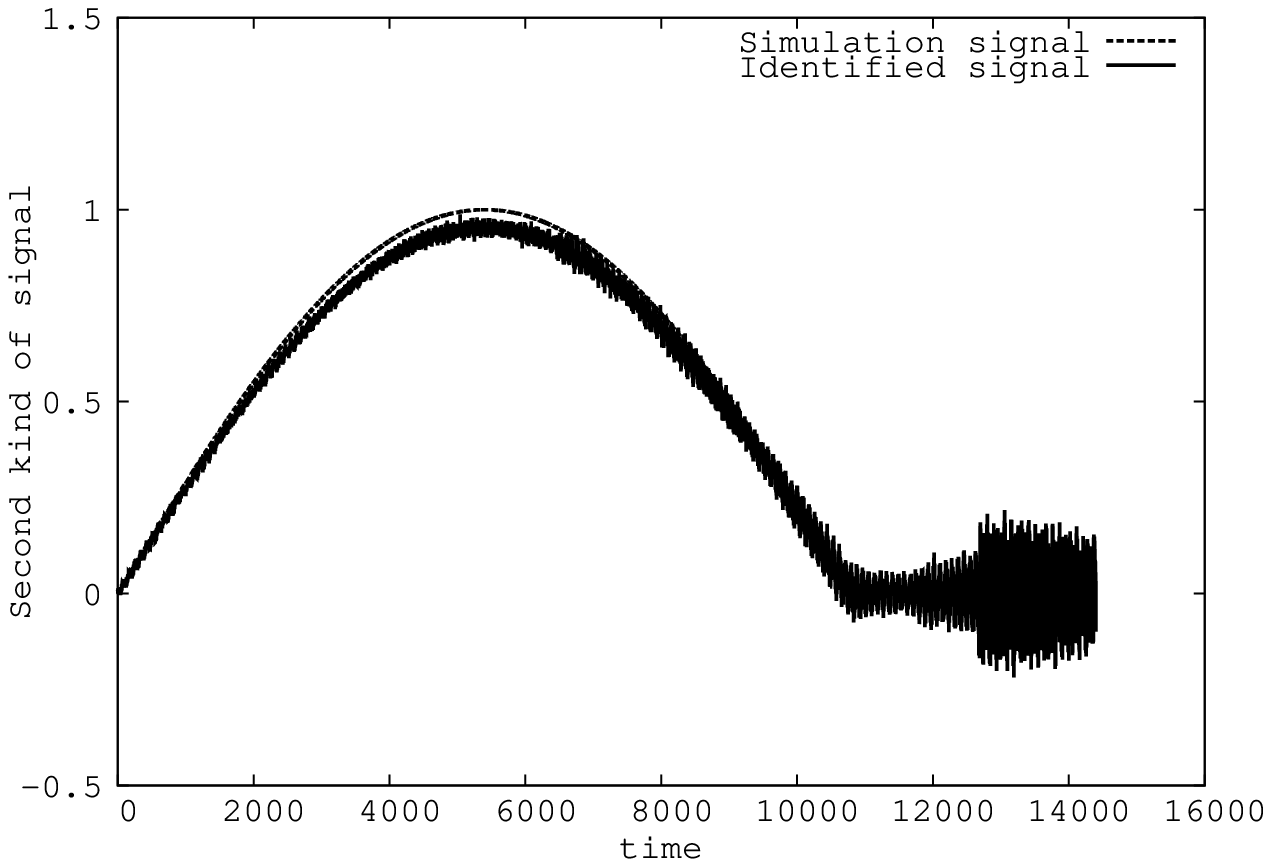}
\end{array}
\end{eqnarray*}
\caption{{\bf{(a)}} Gaussian Noise $0\%$, $\mbox{Error}_2=0,04$ \hspace{0.5cm} {\bf{(b)}} Gaussian Noise $1\%$, $\mbox{Error}_2=0,07$}
\label{Fig3}
\end{center}
\end{figure}

\begin{figure}[H]
\begin{center}
\begin{eqnarray*}
\begin{array}{cccc}
{\bf{(c)}}&&{\bf{(d)}}\\
\includegraphics[width=7cm,angle=0]{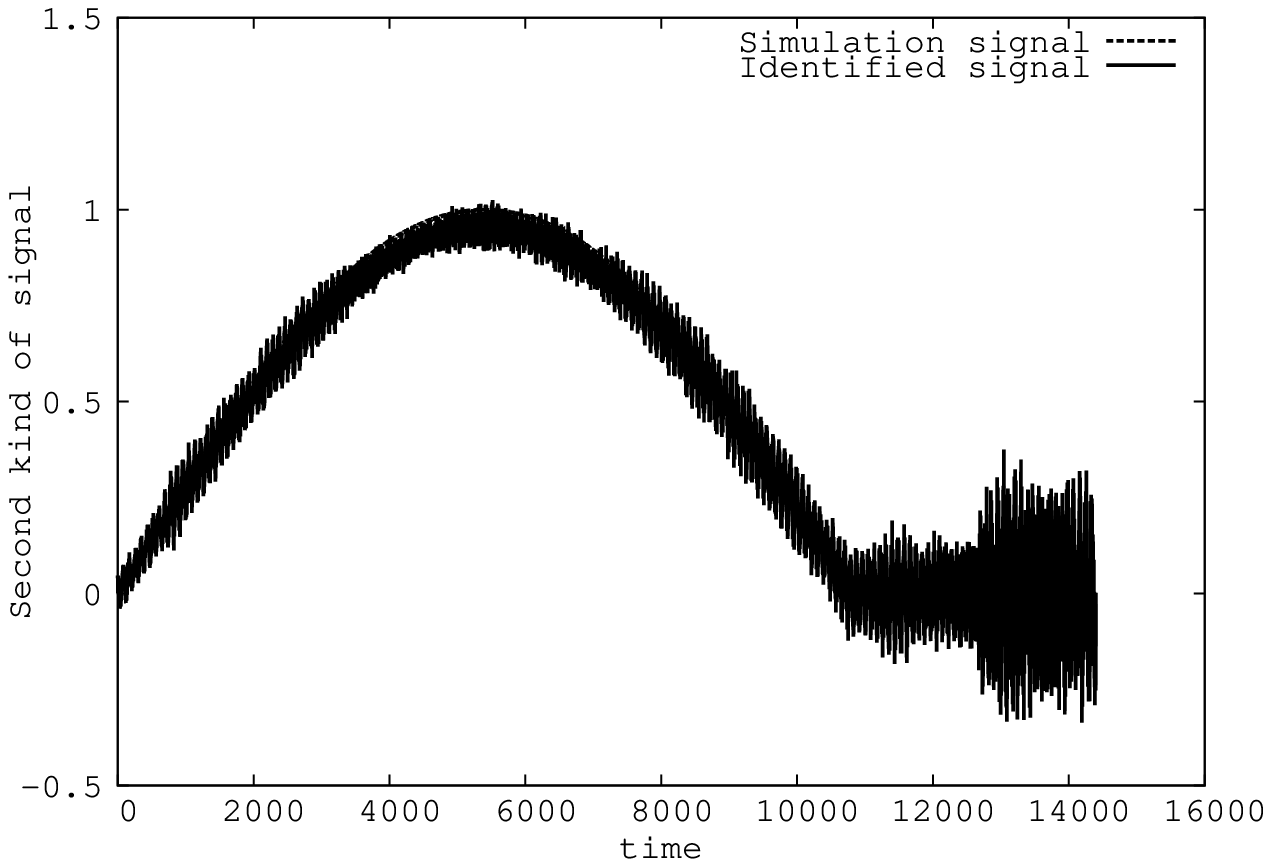}&&\includegraphics[width=7cm,angle=0]{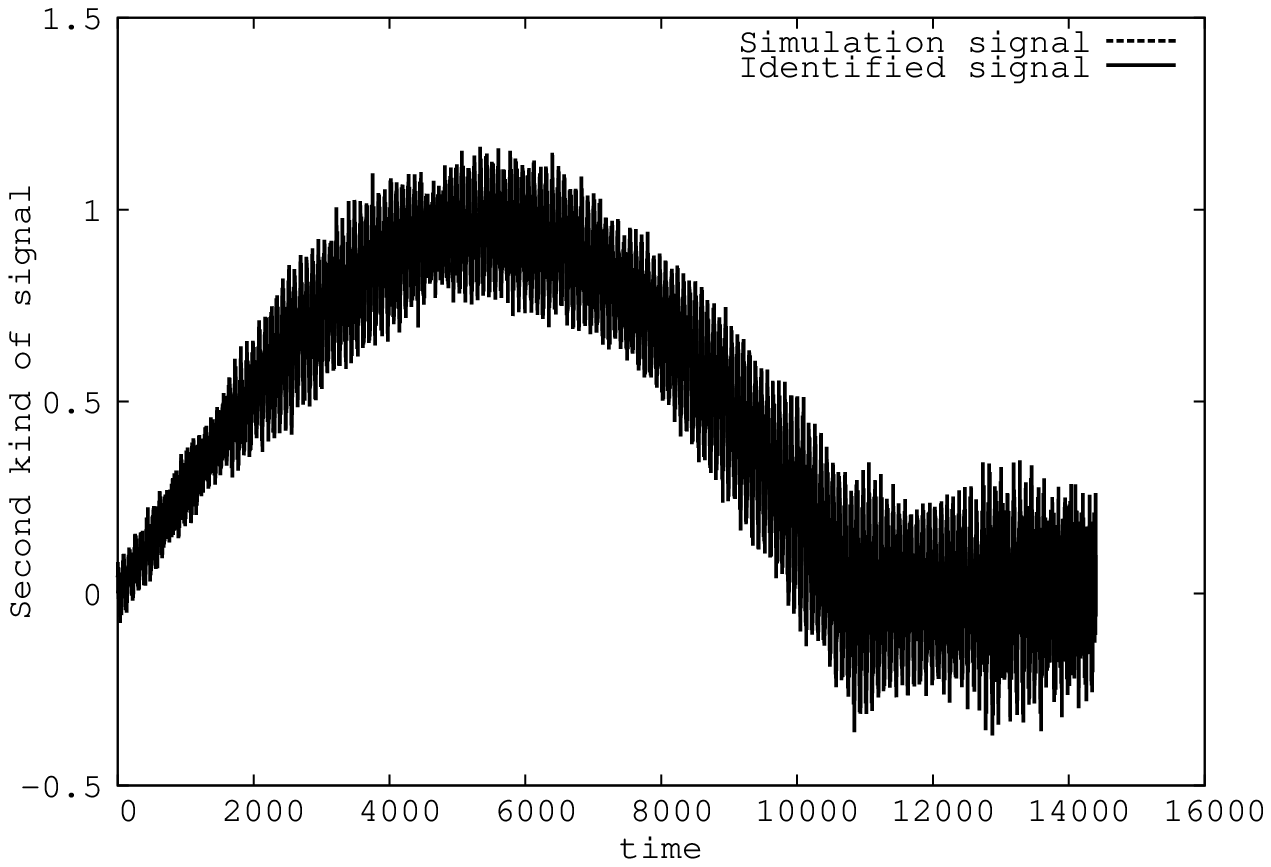}
\end{array}
\end{eqnarray*}
\caption{{\bf{(c)}} Gaussian Noise $3\%$, $\mbox{Error}_2=0,11$ \hspace{0.5cm} {\bf{(d)}} Gaussian Noise $5\%$, $\mbox{Error}_2=0,17$}
\label{Fig4}
\end{center}
\end{figure}

Figures \ref{Fig3} and \ref{Fig4} present the curve, 
labeled "Simulation signal", of the second kind of signal $\lambda_2$ 
defined in (\ref{lambdas}) and the curve of the identified 
signal obtained by applying the deconvolution method on measures, taken 
at the observation node $k=5$, generated by $\lambda=\lambda_2$ and a source 
located at the node $2$ i.e., $S=(0,1,0,0,0)^{\top}$. 

The $L^2$ relative 
error on the identification of $\lambda_2$ while changing the observation 
node $k$ is given by
\begin{table} [H]
\centering
\begin{tabular}{|l|c|c|c|c|c|r|}
   \hline\hline
Observer node     &    $k=1$    &    $k=2$   &   $k=3$    &  $k=4$ & $k=5$ \\ \hline
$\mbox{Error}_2$ & $0,095$ & $0,128$  & $0,111$ & - & $0,109$ \\  
\hline\hline
\end{tabular}
\caption{Identification of $\lambda_2$: Source $S=(0,1,0,0,0)^{\top}$ and Gaussian Noise $3\%$}
\label{tab1}
\end{table}


\begin{figure} [H]
\begin{center}
\begin{eqnarray*}
\begin{array}{cccc}
{\bf{(a)}}&&{\bf{(b)}}\\
\includegraphics[width=7cm,angle=0]{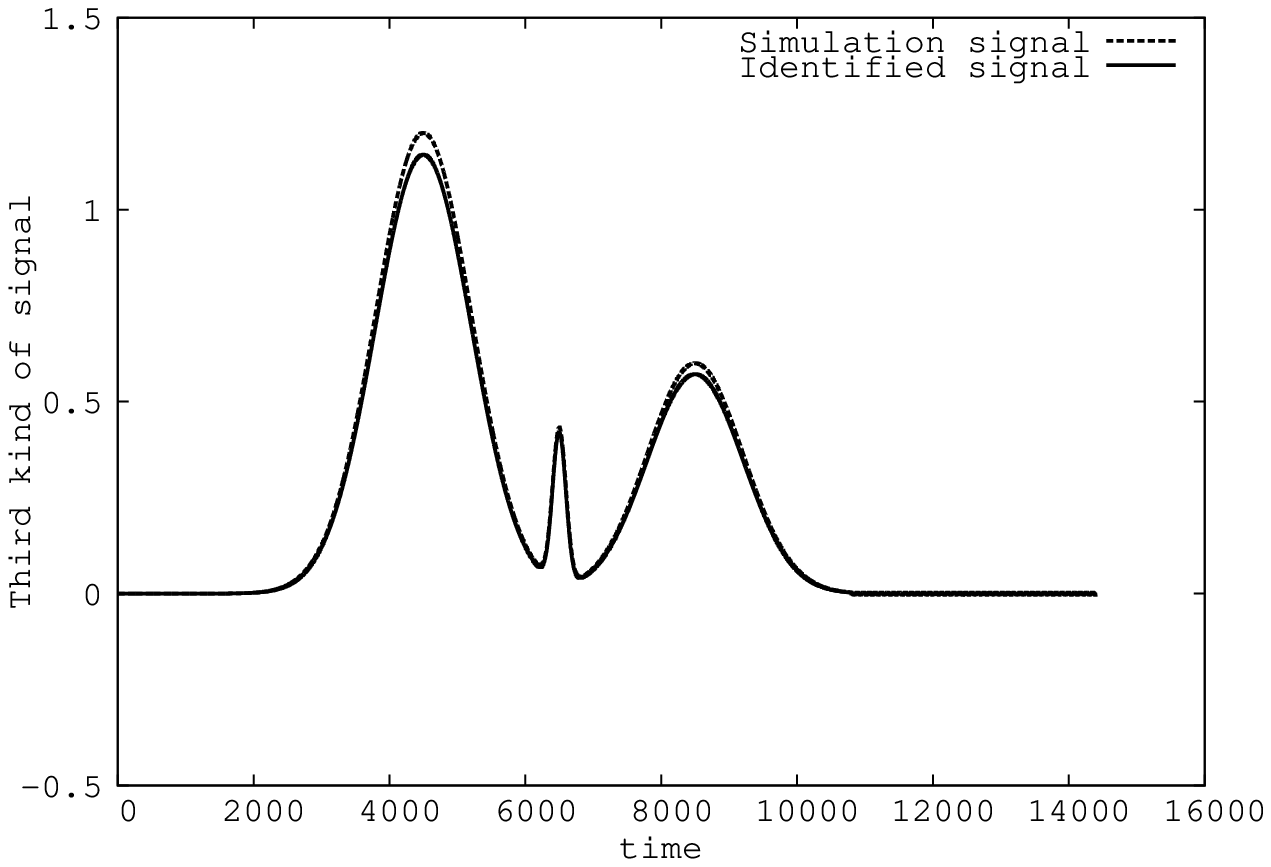}&&\includegraphics[width=7cm,angle=0]{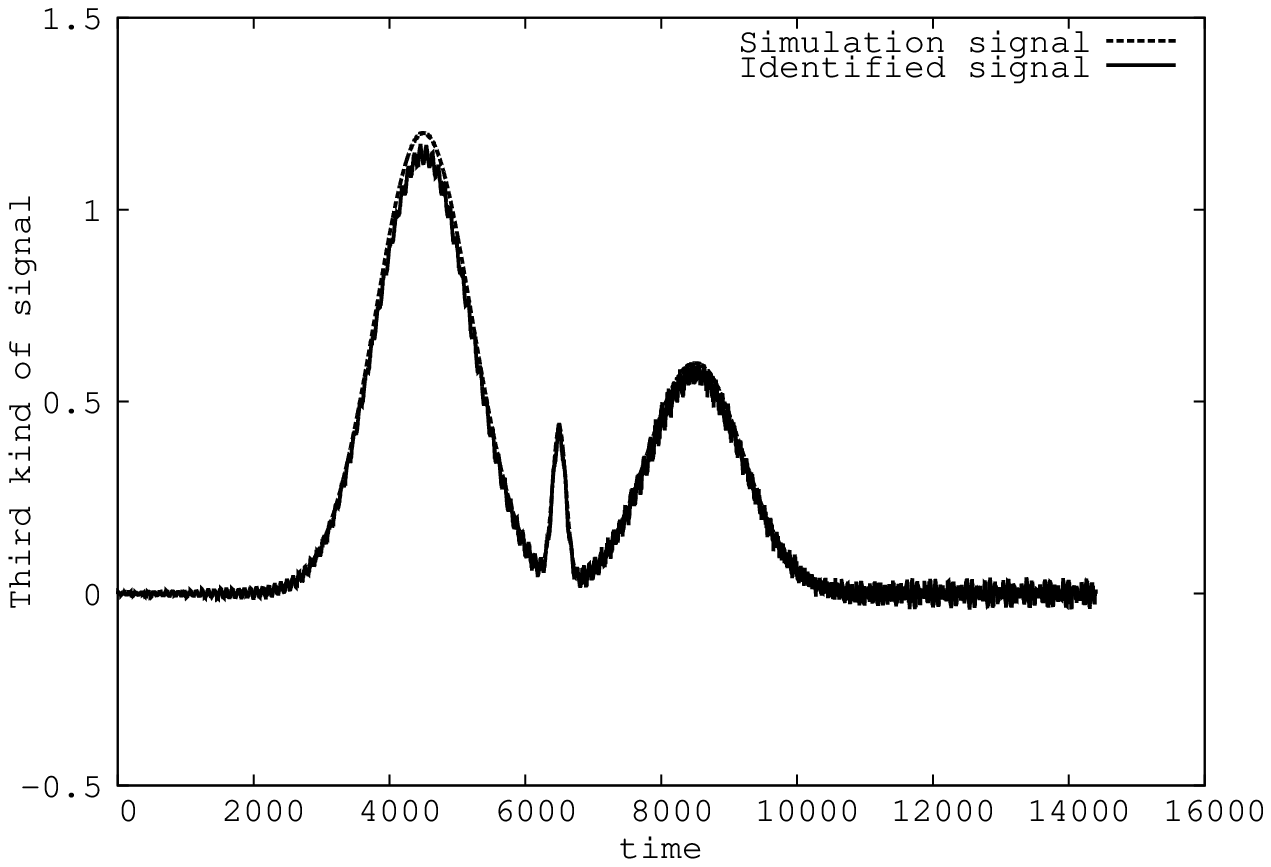}
\end{array}
\end{eqnarray*}
\caption{{\bf{(a)}} Gaussian Noise $0\%$, $\mbox{Error}_3=0,04$ \hspace{0.5cm} {\bf{(b)}} Gaussian Noise $1\%$, $\mbox{Error}_3=0,06$}
\label{Fig5}
\end{center}
\end{figure}

\begin{figure}[H]
\begin{center}
\begin{eqnarray*}
\begin{array}{cccc}
{\bf{(c)}}&&{\bf{(d)}}\\
\includegraphics[width=7cm,angle=0]{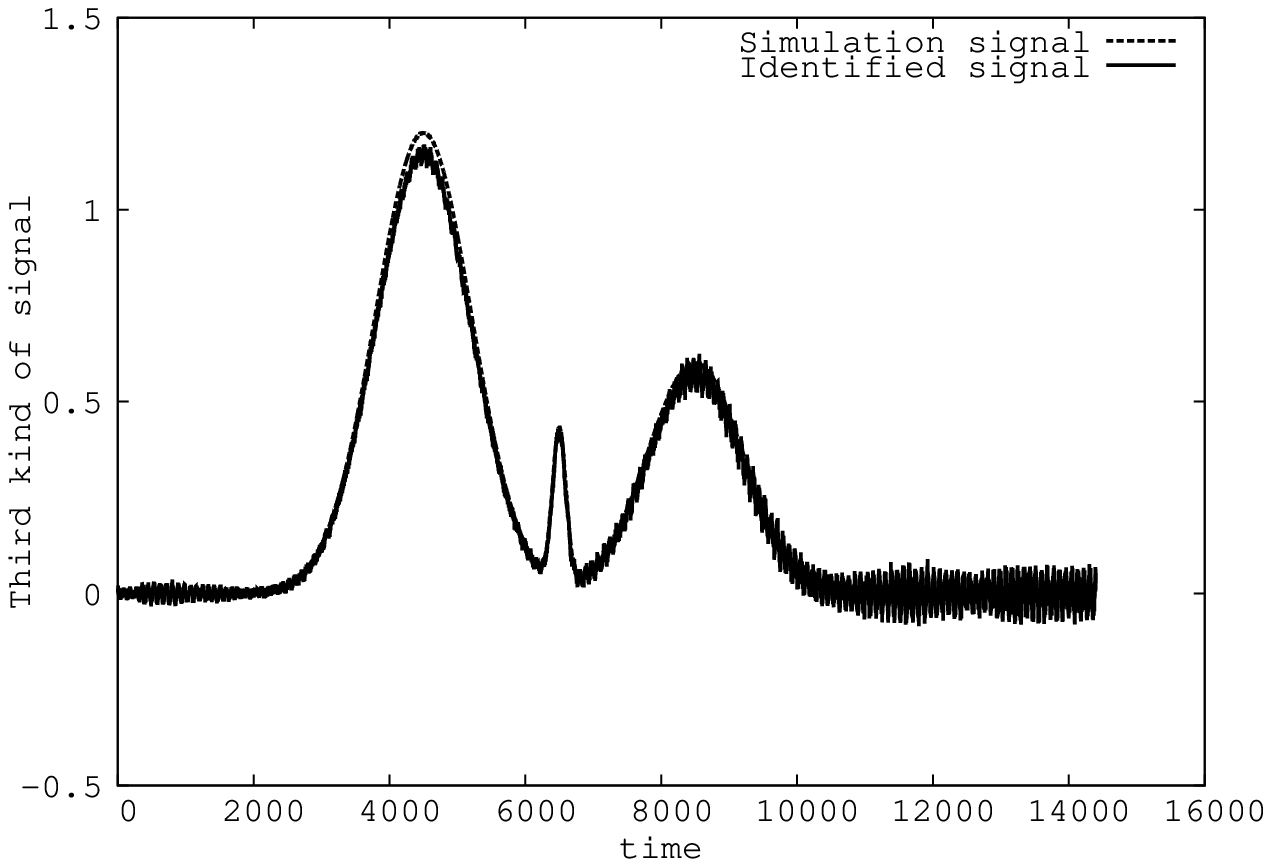}&&\includegraphics[width=7cm,angle=0]{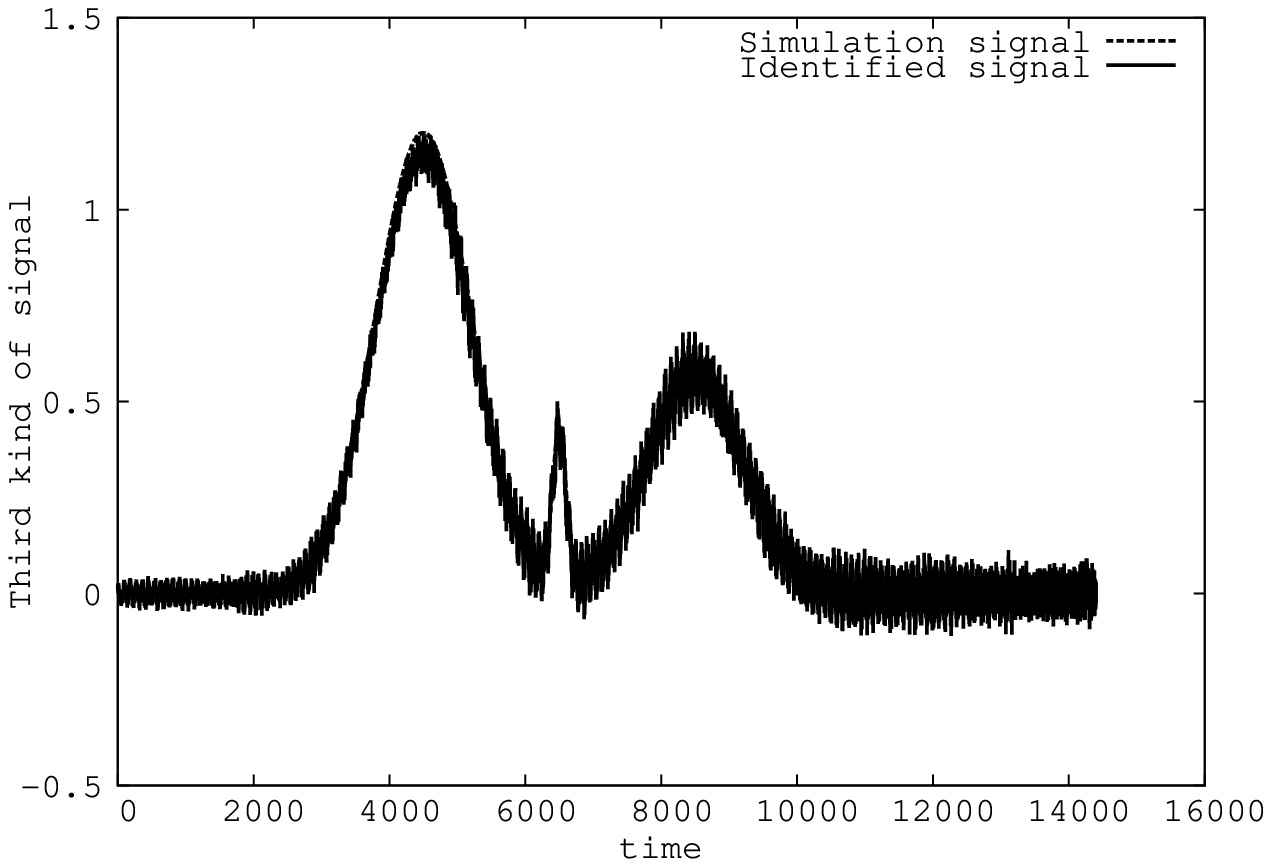}
\end{array}
\end{eqnarray*}
\caption{{\bf{(c)}} Gaussian Noise $3\%$, $\mbox{Error}_3=0,07$ \hspace{0.5cm} {\bf{(d)}} Gaussian Noise $5\%$, $\mbox{Error}_3=0,10$}
\label{Fig6}
\end{center}
\end{figure}
The numerical results in Figures \ref{Fig5} and \ref{Fig6} give the curve, labeled "Simulation signal", of the third kind of signal $\lambda_3$ introduced in (\ref{lambdas}) and the curve of the identified signal obtained by applying the deconvolution method on measures, taken in the observation node $k=3$, generated by $\lambda=\lambda_3$ and a source located at the node $1$ i.e., $S=(1,0,0,0,0)^{\top}$. The $L^2$ relative error on the identification of $\lambda_3$ while changing the observation node $k$ is given by
\begin{table} [H]
\centering
\begin{tabular}{|l|c|c|c|c|c|r|}
   \hline\hline
Observer node     &    $k=1$    &    $k=2$   &   $k=3$    &  $k=4$ & $k=5$ \\ \hline
$\mbox{Error}_3$ & $0,098$ & $0,088$  & $0,071$ & - & $0,113$ \\  
\hline\hline
\end{tabular}
\caption{Identification of $\lambda_3$: Source $S=(1,0,0,0,0)^{\top}$ and Gaussian Noise $3\%$}
\label{tab1}
\end{table}
The numerical results obtained from the identification of the different 
signals $\lambda_{l=1,2,3}$ introduced in (\ref{lambdas}) show that 
the identification method based on deconvolution developed in the present 
study enables to identify the unknown emitted signal $\lambda$ using 
measures taken in any observation node $k$ except the node $k=4$. In fact, 
the numerical difficulties encountered while identifying the emitted signal 
using the observation node $k=4$ were expected since the latter is not 
a strategic node and in view of the application of Titchmarsh's 
Theorem in the last part of the proof of Theorem 4.2.

{The noisy tails observed for added noise are due to the amplification
of the noise by the numerical errors. Despite them, the signal $\lambda(t)$
is well approximated.}

The numerical results presented in this subsection on 
the identification of the emitted signal $\lambda$ were obtained 
solving the triangular linear system (\ref{LS_discrete}).
The diagonal of this system is $\Phi_k(t_1)$ which is close
to $0$. To obtain convergence of the iteration, 
we replaced the term 
$\Phi_k(t_1)$ by $\Phi_k(t_1) + r $ with $r>0$ which 
can be seen as a regularisation parameter. In the absence of noise,
$r \approx 1$ while with the noise present we took
$r \approx 10$.

At the beginning of the subsection, we indicate that $\lambda(t)$
was chosen to vary slowly compared to the typical periods of the
network. This is the most favorable situation to reconstruct
$\lambda(t)$. When this forcing varies on time scales comparable
to the typical periods, the identification of $\lambda(t)$ can 
also be done but it is less accurate.

\subsubsection{Identification of $\lambda$ using Fourier basis}

As introduced in section 4, the second way of identifying the source 
signal $\lambda$ consists in using the over-determined linear system 
(\ref{adj2})-(\ref{Pm_ij}) once the source position $S$ is localized 
to compute the coefficients 
$\lambda_m=\langle \lambda,\varphi_m\rangle_{L^2(0,T)}$ for $m=1,\dots,M$. 
To this end, since the source position was determined in the previous 
subsection i.e., $S=\big(0,0,1,0,0\big)^{\top}$ it follows that in 
the linear system (\ref{id_sys})-(\ref{Pm_ij}) the third equation 
contains the source term $\lambda_m$. Therefore, we proceed as follows:

For $m=1$ to $M$ do
\begin{itemize}
\item Compute the right-hand side vector $P_m^{1,2}$ in (\ref{Pm_ij}).

\item Solve $3$ equations other than the third equation (which contains the 
source term) from the over-determined linear 
system (\ref{id_sys}) and compute $\bar X_m^{1,2}$.

\item Use $\bar X_m^{1,2}$ in the third equation of the system (\ref{id_sys}) to deduce the coefficient $\lambda_m$.

\item End For
\end{itemize}
Hence, the identified source signal is defined by 
$$\lambda(t)=\sum_{m=1}^M \lambda_m \varphi_m(t),$$ 
for all $t\in(0,T)$. 
We applied this procedure to the data generated in section 6.2 and 
obtained an approximation of the signal, using the first seven harmonics.
For this rapid $\lambda(t)$, both the Fourier and the deconvolution 
methods provided an estimate of the signal.

In the numerical experiments that we carried out, the Fourier 
method of reconstructing the unknown source intensity function $\lambda$ 
appeared to work well. Nevertheless, this second method admits the two 
following disadvantages with respect to the deconvolution method. 
First, to obtain an accurate approximation of $\lambda$, we need a
relatively large number $M$ of eigenfunctions $\varphi_m$. This involves 
functions of high frequency which in turn affect the accuracy of the method.
Another disadvantage is that to determine each component 
$\lambda_m=\langle \lambda,~
\varphi_m\rangle_{L^2(0,T)}$, we need to solve a linear system of 
$(N-2)$ equations. This increases the total reconstruction cost of 
$\lambda$ using this second method.

On the other hand, we noted some convergence problems for the 
deconvolution method on both the fast signal of section 6.2 and
the slow signal of section 6.3.1. The linear system 
\eqref{Lin_Syst_lam_n} needs 
to be regularized to obtain convergence. The Fourier method appears
to be more fail-proof. A systematic study is needed to
determine -depending on the typical periods of oscillation of
the network and the period of $\lambda(t)$- what reconstruction
method will be more appropriate. Both methods should be used in
practical cases.
 

\section{Conclusion }

We studied the nonlinear inverse source problem of localizing a 
node source emitting an unknown time-dependent signal in a network 
defined over a $N-$ node graph, on which propagates a miscible flow.

We proved under reasonable assumptions that time records of 
the associated state taken in a strategic set of two nodes yield 
uniqueness of the two unknown elements defining the occurring node 
source. Graph theory arguments guide us in the choice of observation nodes
for a given network.

The article also presents a direct identification method that localizes 
the position of the node source in the network by solving a set of 
well posed linear systems and proposes two different algorithms 
to identify its unknown emitted signal.

Numerical experiments on a $5$-node graph showed
the effectiveness of the identifications of the node source and the
signal.

Due the wide spectrum of applications covered by the 
inverse source problem studied in this article, perspectives are numerous.
A challenge will be to find an efficient numerical procedure to
identify sources on networks of large size. 

\section*{Acknowledgment}

This work is part of the XTerM project, co-financed by the European 
Union with the European regional development fund (ERDF) and by 
the Normandie Regional Council. We also acknowledge partial support from
the ANR grant "Fractal grid".



{\bf{Appendix}}\\

In the sequel, we establih the proof of Lemma \ref{Null_State}.\\

{\bf{Proof.}} From expanding the solution $X$ of the system (\ref{Lemma_1}) in the orthonormal family $\{v^1,\dots,v^N\}$ of $\mathbb R^N$ and using the notations: $y_n(T)=\langle X(T),v^n\rangle$ and $\dot y_n(T)=\langle \dot X(T),v^n\rangle$ for $n=1,\dots,N$, we obtain for all $t\in (T^*,T)$

\begin{eqnarray}
X(t)=\Big(y_1(T) + (t-T)\dot y_1(T)\Big) v^1 + \displaystyle \sum_{n=2}^{N} \Big(y_n(T)\cos\big(\omega_n(t-T)\big)  + \displaystyle \frac{\dot y_n(T)}{\omega_n}\sin\big(\omega_n(t-T)\big)\Big)v^n
\label{Express_X}
\end{eqnarray}

Since, we have $x_{k}(t)=0$ for all $t\in(T^*,T)$ and $k\in\{k_1,\dots,k_{\ell}\}$ it follows in view of (\ref{Express_X}) and by analytic continuation that for all $t\in]T^*,+\infty[$ and all $k\in\{k_1,\dots,k_{\ell}\}$,

\begin{eqnarray}
\Big(y_1(T) + (t-T)\dot y_1(T)\Big) v^1_{k} + \displaystyle \sum_{n=2}^{N} \Big(y_n(T)\cos\big(\omega_n(t-T)\big) + \displaystyle \frac{\dot y_n(T)}{\omega_n}\sin\big(\omega_n(t-T)\big)\Big)v^n_{k}=0
\label{X_strat}
\end{eqnarray}

Let $\tau>0$. By integrating the equation (\ref{X_strat}) over $(T,T+\tau)$, we find 

\begin{eqnarray}
\Big(\tau y_1(T) + \frac{\tau^2}{2}\dot y_1(T)\Big) v_{k}^1 + \displaystyle \sum_{n=2}^{N} \Big(y_n(T) \frac{\sin\big(\omega_n\tau\big)}{\omega_n} + \displaystyle \frac{\dot y_n(T)}{\omega_n}\frac{1 - \cos\big(\omega_n\tau\big)}{\omega_n}\Big)v_{k}^n=0, \; \forall \tau>0
\label{int_X_strat}
\end{eqnarray} 

Then, from dividing by $\tau^2$ the left hand side in (\ref{int_X_strat}) and setting the limit when $\tau$ tends to $+\infty$, we obtain: $\dot y_1(T)v_{k}^1=0$ for all $k\in\{k_1,\dots,k_{\ell}\}$. Furthermore, using this last result in (\ref{int_X_strat}) and dividing now by $\tau$ with setting the limit when $\tau$ tends to $+\infty$ gives: $y_1(T)v_{k}^1=0$ for all $k\in\{k_1,\dots,k_{\ell}\}$. Thus, in (\ref{X_strat}) the term involving $v^1_{k}$ vanishes. Let $n_0\in\{2,\dots,N\}$. From multiplying (\ref{X_strat}) by $\mathrm{cos}\big(\omega_{n_0}(t-T)\big)$ and integrating over $(T,T+\tau)$ then, dividing the obtained result by $\tau$ and setting the limit when $\tau$ tends to $+\infty$ we get: $y_{n_0}(T) v_{k}^{n_0}=0$ for all $k\in\{k_1,\dots,k_{\ell}\}$. Therefore, we find for all $n\in\{1,\dots,N\}$

\begin{eqnarray}
\begin{array}{lll}
\qquad y_n(T) v_{k}^n=0, \;\;  \forall \; k\in\{k_1,\dots,k_{\ell}\} \quad\implies\quad y_1(T)=\dots=y_N(T)=0\\
\qquad \dot y_1(T) v_{k}^1=0,\;\;  \forall \; k\in\{k_1,\dots,k_{\ell}\} \quad\implies\quad \dot y_1(T)=0
\end{array}
\label{X_initial}
\end{eqnarray}

The two implications in (\ref{X_initial}) are obtained since the set $\{k_1,\dots,k_{\ell}\}$ is {\it{strategic}}. In addition, using (\ref{X_initial}) in (\ref{Express_X}) and deriving the obtained form with respect to $t$ gives

\begin{eqnarray}
\dot X(t)=\displaystyle \sum_{n=2}^{N} \dot y_n(T) \cos\big(\omega_n(t-T)\big)v^n,\quad \forall t\in(T^*,T)
\label{Express_dot_X}
\end{eqnarray}

Moreover, as we have: $\dot x_{k}(t)=0$ for all $t\in (T^*,T)$ and $k\in\{k_1,\dots,k_{\ell}\}$ then, by analytic continuation it follows from (\ref{Express_dot_X}) that for all $k\in\{k_1,\dots,k_{\ell}\}$

\begin{eqnarray}
\displaystyle \sum_{n=2}^{N} \dot y_n(T) \cos\big(\omega_n(t-T)\big)v_{k}^n=0, \qquad \forall t\in ]T^*,+\infty[
\label{Express_dot_X_kl}
\end{eqnarray}  

Let $n_0\in\{2,\dots,N\}$ and $\tau>0$. From multiplying (\ref{Express_dot_X_kl}) by $\mathrm{cos}\big(\omega_{n_0}(t-T)\big)$ and integrating over $(T,T+\tau)$  then, dividing the obtained result by $\tau$ and setting the limit when $\tau$ tends to $+\infty$ we obtain $\dot y_{n_0}(T) v_{k}^{n_0}=0$ for all $k\in\{k_1,\dots,k_{\ell}\}$. Hence, we have for all $n\in\{2,\dots,N\}$:

\begin{eqnarray}
\dot y_n(T) v_{k}^n=0, \;\; \forall \; k\in\{k_1,\dots,k_{\ell}\}  \quad\implies\quad \dot y_2(T)=\dots=\dot y_N(T)=0
\label{dot_X_initial}
\end{eqnarray}

The implication in (\ref{dot_X_initial}) is obtained since $\{k_1,\dots,k_{\ell}\}$ is {\it{strategic}}. Therefore, (\ref{X_initial}) and (\ref{dot_X_initial}) imply that $X(T)=\dot X(T)=\vec 0$ which is the result announced in (\ref{Result_lemma}). \hspace{2cm}$\blacksquare$

\end{document}